\documentclass[11pt]{amsart}
\usepackage{a4wide}
\usepackage{amssymb}
\usepackage{amsthm}
\usepackage{amsmath}
\usepackage{amsopn}
\usepackage{array}
\usepackage{verbatim}

\DeclareMathOperator{\supp}{supp}

     \newcommand{\NN}{\mathbb{N}}

     \newcommand{\ZZ}{\mathbb{Z}}

\theoremstyle{plain}
\newtheorem{teo}{Teorem}
\newtheorem{theorem}[teo]{Theorem}
\newtheorem{proposition}[teo]{Proposition}
\newtheorem{lemma}[teo]{Lemma}
\newtheorem*{tools}{Fundamental Tools}

\theoremstyle{definition}
\newtheorem{definition}[teo]{Definition}

\theoremstyle{remark}
\newtheorem{remark}[teo]{Remark}

\newcommand{\Hm}[1]{\leavevmode{\marginpar{\tiny%
$\hbox to 0mm{\hspace*{-0.5mm}$\leftarrow$\hss}%
\vcenter{\vrule depth 0.1mm height 0.1mm width \the\marginparwidth}%
\hbox to
0mm{\hss$\rightarrow$\hspace*{-0.5mm}}$\\\relax\raggedright #1}}}

\begin{document}

\title[Sharpness of the phase transition and exponential decay of cluster size]{Sharpness of the phase transition and exponential decay of the subcritical cluster size for percolation on quasi-transitive graphs}
\author[T.~Antunovi\'c]{Ton\'ci Antunovi\'c}

\author[I.~Veseli\'c]{Ivan Veseli\'c}
\address{Emmy-Noether-Programme of the Deutsche Forschungsgemeinschaft\vspace*{-0.3cm} }
\address{\& Fakult\"at f\"ur Mathematik,\, 09107\, TU\, Chemnitz, Germany   }
\urladdr{www.tu-chemnitz.de/mathematik/schroedinger/members.php}

\curraddr[T.~Antunovi\'c]{Department of Mathematics, UC Berkeley, Berkeley CA 94720, USA}
\curraddr[I.~Veseli\'c]{Institut f\"ur Angewandte Mathematik, 53115 Universit\"at Bonn, Germany}


\begin{abstract}
We study homogeneous, independent percolation on general quasi-transitive graphs.
We prove that in the disorder regime where all clusters are finite almost surely,
in fact the expectation of the cluster size is finite. This extends a well-known theorem
by Menshikov and Aizenman \& Barsky to all quasi-transitive graphs. Moreover we deduce that in this disorder regime
the cluster size distribution decays exponentially, extending a result of  Aizenman \& Newman.
Our results apply to both edge and site percolation, as well as long range (edge) percolation.
The proof is based on a modification of the Aizenman \& Barsky method.
\end{abstract}
\thanks{\jobname.tex; \today }
\maketitle

\section{Introduction}
Percolation theory is devoted to the study of geometric properties of random subgraphs of a given graph.
In particular, one wants to understand which disorder regimes exhibit the existence of an infinite cluster,
i.e.~an infinite component of the subgraph generated by the percolation process.

For percolation models on graphs the low density phase is often defined as the regime
of randomness where the expected cluster-size  is finite, whereas the high density phase
is defined as the disorder regime where there exists an infinite cluster almost surely.
More specifically, for identically distributed, independent models 
there is only one scalar disorder parameter (usually denoted by $p$) which measures the extent
of the randomness. If one denotes the supremum of the parameter values which correspond to
the low density regime by $p_T$, and the infimum of the parameter values which correspond to the high
density regime by $p_H$, then the statement $p_T = p_H$ is called \emph{sharpness of the phase transition}.
In other words, an intermediate phase between the low and high density regime reduces to (at most) a single value
of the parameter $p$.

This result has been proven by Menshikov in \cite{Menshikov-86}
(see also \cite{MenshikovMS-86} by Menshikov, Molchanov \& Sidorenko)
and Aizenman \& Barsky in \cite{AizenmanB-87} for a large class of percolation processes on graphs.
More precisely the results of \cite{Menshikov-86,MenshikovMS-86} cover independent
site percolation on quasi-transitive graphs of subexponential growth. The percolation parameter can 
be different for the different classes of vertices.
In \cite{MenshikovMS-86}, where the proof of \cite{Menshikov-86} is explained in more detail,
it is noted in Remark 6.1 that  the method of proof works also with a relaxed growth condition on the graph,
however, that it is not possible to eliminate it completely.
The results of \cite{AizenmanB-87} hold for directed and undirected, independent, site and bond, short range 
percolation models on  $\ZZ^d$. 
Since the considered graphs are essentially Cayley graphs of $\ZZ^d$ 
their volume growth is polynomially bounded.
The results of \cite{AizenmanB-87} furthermore apply to so called \emph{long range}  edge percolation,
a model where edges may be present between any pair of vertices, with probability decreasing 
in the distance between the vertices. 

Any edge percolation process can be transformed into a site percolation process
by passing to the line graph. If the original graph has a finite number of 
edge orbits under the automorphism group action, the resulting line graph will be quasi-transitive.
Thus the results of \cite{Menshikov-86,MenshikovMS-86} apply to edge percolation, too.
In contrast to this, if we transform a long range edge percolation process to a site percolation process
via  the line graph construction we lose  quasi-transitivity. More precisely, to avoid triviality for
long range edge percolation we need to have an infinite number of different percolation parameters assigned to the 
edges. Thus the edge set decomposes into an infinite number of classes, which means that the 
line graph will have infinitely many vertex types, violating the quasi-transitivity property.

We adapt the method of differential inequalities
used in \cite{AizenmanB-87} for the study of $\ZZ^d$-models, to show that the sharpness of the phase transition actually 
holds for all quasi-transitive graphs. Again we can treat short range site and edge percolation, 
as well as long range edge percolation. On the technical level the differences to \cite{AizenmanB-87} 
are the following: In \cite{AizenmanB-87} finite torus graphs are used to approximate the infinite $\ZZ^d$ graph,
which has the advantage that the approximating  graph is still homogeneous, i.e. transitive.
In the general case this is not possible, thus we work with finite approximation graphs which have a `boundary'.
As a consequence of this and quasi-transitivity rather than transitivity, in comparison to \cite{AizenmanB-87} additional finite volume terms appear 
in the differential inequalities. We control these correction terms to show that the 
modified inequalities still lead to a proof of the sharpness of the phase transition.

All graphs mentioned so far have a rich algebraic structure which is formulated in terms of transitivity
or quasi-transitivity. 
Looking at the modifications of the Aizenman \& Barsky method needed for quasi-transitive graphs
one gets the impression that one can adopt the proof to obtain the same results for 
graphs which have uniform local combinatorial complexity bounds,  
without having necessarily a large automorphism group. An example would be the Penrose tiling, 
whose percolation properties were studied by Hof in \cite{Hof-98b}. 
For such graphs the finite volume effects can probably be controlled
in a similar way as for quasi-transitive graphs. M\"uller and Richard have recently in \cite{MuellerR} 
obtained related results 	using different techniques.

Apart from the fact that the equality $p_T = p_H$ establishes for homogeneous, independent models
on $\ZZ^d$ that percolation has only one phase transition,
it has also played a crucial role for the proof of Kesten's theorem \cite{Kesten-80}, namely,  that on the
two dimensional lattice $\ZZ^2$ we have $p_T=p_H = 1/2$ for the edge percolation model.

\medskip

Closely related to the sharpness of the phase transition are the results \cite{AizenmanN-84}
by Aizenman \& Newman. On the one hand they prove the divergence of the expectation value of the cluster size
as the disorder parameter approaches $p_T$ from below, a statement which is used in \cite{AizenmanB-87} for the proof of
$p_T=p_H$. On the other hand, Aizenman \& Newman show the exponential decay of the cluster size distribution for
$p<p_T$. We deduce that these results hold actually for all quasi-transitive graphs.

%
%

\medskip

Our interest in the sharpness of the phase transition stems from the study of
percolation Hamiltonians, more precisely adjacency and combinatorial Laplacians on percolation subgraphs,
and in particular their integrated density of states (IDS).
While for the definition of the IDS for general graphs with a free, amenable, quasi-transitive group action,
cf.~\cite{Veselic-05a,Veselic-05b}, an understanding  of the phase transition(s) is not necessary,
it seems that for the proof of Lifshitz asymptotics such understanding is crucial.

Lifshitz tails describe the  asymptotic behaviour of the IDS near the boundaries of the spectrum
and have been established for independent site percolation on $\ZZ^d$ by Biskup \& K\"onig \cite{BiskupK-01a} 
and independent edge percolation on $\ZZ^d$ by Kirsch \& M\"uller, resp.~M\"uller \& Stollmann in \cite{KirschM-06,MuellerS}. Before that Klopp and Nakamura \cite{KloppN-03} derived partial results for the random hopping model, which includes the edge percolation Hamiltonian as a special case.
In \cite{AntunovicV-b,AntunovicV-c} we study subcritical independent site and edge percolation on amenable Cayley graphs.
For the adjacency and the combinatorial Laplace operator we obtain the asymptotics of the IDS at the spectral edge.
Our results depend on the decay rate of the cluster size distribution.
They apply in particular to Lamplighter graphs which are amenable, but have exponential volume growth.
Let us also remark, that the results of \cite{AntunovicV-c} cover combinatorial Laplacians
on long range percolation graphs in the subcritical phase. This explains our objective to derive the exponential decay of the subcritical cluster size 
also for graphs with arbitrary growth behaviour, and for long range edge percolation.

\bigskip

The structure of this paper is as follows: In the next section we define our percolation model
and state the main results. In Section \ref{s-basics} we state some basic facts which are common for site
and for edge percolation. 
The two subsequent sections are devoted to (short range) edge percolation. 
Namely, in Section \ref{s-order-parameter} we introduce the (finite volume) order parameters and in Section \ref{s-diffinequalities}
we show that they obey certain differential inequalities. In Section \ref{s-site-model} we establish the same facts for site percolation.
This allows us to complete the proof of our main result for both types of percolation processes in Section \ref{s-completion}.
The last section contains an extension of our main result to long range edge and oriented percolation.

\section{Notation and results}
Let $G=(V,E)$ be an infinite, countable,
connected graph, with vertex set $V$ and edge set $E$.
The fact that the vertices $x$ and $y$ are adjacent will be
denoted by $x \sim y$ and $[x,y]$ will stand for the unoriented
edge which connects $x$ and $y$. By $d \colon V\times V \to
\mathbb{R}$ we denote the usual graph distance, that is $d(x,y)$
is the length of a shortest path between two vertices $x$ and $y$.
For a vertex $x$ and a nonnegative integer $n$, $B(x,n)$ is the
ball, with center $x$, of radius $n$, in the above metric. For the
sphere of radius $n$ around $x$  we shall write $S(x,n):=\left\{y
\in V|d(y,x)=n\right\}$. The group of graph automorphisms will be
denoted by $Aut(G)$. For any subgraph $G'$, $|G'|$ will stand for
the number of vertices in $G'$, which may be infinite. When a
subgraph $G'$ is given, we will say that two neighboring vertices
of $G$, $x$ and $y$, are \emph{directly connected in} $G'$, if the
edge $[x,y]$ is an edge of the graph $G'$. For two subgraphs
$G_{1} = (V_{1},E_{1})$ and $G_{2} = (V_{2},E_{2})$ we define
their intersection $G_{1} \cap G_{2} := (V_{1} \cap V_{2}, E_{1}
\cap E_{2})$. Notice that the intersection is always well defined.
The notation $G' \subset G$ means that $G'$ is either a proper
subgraph of $G$ or $G$ itself.

A graph $G$ is called \emph{quasi-transitive}, if there exists a
finite set of vertices $\mathcal{F}$ such that for any vertex $x$
there is a $y \in \mathcal{F}$ and $\gamma \in Aut(G)$ such that
$\gamma y =x$. In the following we will always assume that the set
$\mathcal{F}$ is minimal with respect to inclusion. Such an
$\mathcal{F}$ will be called \emph{fundamental domain}. For any
graph $G$ we can consider the action of the group $Aut(G)$ on the
set of vertices $V$. Thus, a graph is quasi-transitive, if and only
if the set of vertices is decomposed into finitely many orbits,
with respect to this action. A fundamental domain is then any set
of vertices which intersects each orbit in exactly one vertex.
The number of elements in any fundamental domain is the same. If a
fundamental domain contains only one element we call the graph
\emph{transitive}.

Now we introduce the usual nearest neighbor Bernoulli bond
percolation model. We fix some parameter $0 \leq p \leq 1$. For
each edge of the graph $G$ we say that it is open with probability
$p$ and closed with probability $1-p$, independently of all other
edges. That is, for each edge $e \in E$ we take a probability
space $(\Omega_{e}, P(\Omega_{e}), \mathbb{P}_{e})$, where
$\Omega_{e}=\left\{0,1\right\}$, $P(\Omega_{e})$ is the power set
of $\Omega_{e}$ and $\mathbb{P}_{e}(1)=p$,
$\mathbb{P}_{e}(0)=1-p$. The \emph{percolation probability space}
$(\Omega,\mathcal{A},\mathbb{P})$ is defined as the product of
these probability spaces, that is $\displaystyle \Omega:=\prod_{e
\in E}\Omega_{e}$, $\displaystyle \mathcal{A}:=\bigotimes_{e \in
E}P(\Omega_{e})$ and $\displaystyle \mathbb{P}:=\bigotimes_{e \in
E}\mathbb{P}_{e}$. The probability measure $\mathbb{P}$ obviously
depends on $p$. This dependence will sometimes be stressed by
writing $\mathbb{P}_{p}$ instead of $\mathbb{P}$. The same holds
for the expectation $\mathbb{E}$.

Elements of $\Omega$ will be called \emph{configurations} because
each of them uniquely represents some configuration of open and
closed edges. For a given configuration $\omega$ and a given edge
$e$, the value $\omega_{e}$ will be called the \emph{state of
$e$}. By $G(\omega)$ denote the subgraph of $G$ obtained by
deleting all closed edges with respect to the configuration
$\omega$, i.e.~for the set of vertices of $G(\omega)$ we take the
set of vertices of the graph $G$, while the set of edges of
$G(\omega)$ is the set of open edges with respect to the
configuration $\omega$. Connected components of $G(\omega)$ are
called \emph{clusters}. The cluster containing the vertex $x$ will
be denoted by $C_{x}(\omega)$. The probability measure is
invariant under the graph automorphisms and so, in the case of a
transitive graph, the probabilistic properties of $C_{x}(\omega)$
do not depend on the choice of $x$. Thus, in this case, we
will often assume that a certain vertex $x$ is fixed and denote
the cluster $C_{x}(\omega)$ by $C(\omega)$.

The nearest neighbor site percolation model is introduced in an analogous way. For each vertex $x$ we say it is open with some probability $p \in [0,1]$ and otherwise closed, independently of all other vertices. In other words, we consider the probability space $(\Omega, \mathcal{A}, \mathbb{P}) := \prod_{x \in V}(\Omega_{x}, P(\Omega_{x}), \mathbb{P}_{x})$, where $(\Omega_{x}, P(\Omega_{x}), \mathbb{P}_{x})$ is defined in the same way as $(\Omega_{e}, P(\Omega_{e}), \mathbb{P}_{e})$ before. Now for some given configuration $\omega$ the percolation graph $G(\omega)$ is defined simply as the subgraph induced by the set of open vertices with respect to the configuration $\omega$. Clusters are again defined as connected components of $G(\omega)$. We will use the same notation as in the bond model. Note that subgraphs $G(\omega)$ do not have to contain all the vertices of $G$ and thus the event $\left\{|C_{x}|=0\right\}$ has positive probability (namely equal to $1-p$).

\begin{remark}\label{rem: 1}
Everything we mentioned above can be defined for any subgraph
$G'=(V',E')$. Of course, the probability space will be different, but
every event $T$ in the new probability space can (and will) be
identified with the cylinder set $T \times \prod_{e \notin
E'}\Omega_{e}$ ($T \times \prod_{x \notin V'}\Omega_{x}$ in the case of the site model). In the corresponding probability spaces these
events have the same probabilities. Thus we use the same notation
for the corresponding probability measures. Since the notion of
clusters in $G'$ and $G$ is not the same, we will denote by
$C_{x}^{G'}$ the cluster of $x$ in the graph $G'$.
\end{remark}

Since the statements in the present section hold equally for site and for bond percolation, we will use in this section simply the term percolation.
Next we will describe the most basic properties of the percolation process, without giving proofs. In Remark \ref{rem: vertex independence} in the next section we briefly
sketch how these properties are proven.

An important property of  percolation is the existence of
a \emph{phase transition} between a percolating and a non-percolating phase.
Consider some fixed vertex $x$ and the
event $\left\{|C_{x}|=\infty\right\}$. The probability of this
event $\mathbb{P}_{p}(|C_{x}|=\infty)$ is equal to $0$ when $p=0$
and $1$ when $p=1$. Furthermore, it can be shown that
$\mathbb{P}_{p}(|C_{x}|=\infty)$ is a non-decreasing function of
$p$. Thus, if we define $p_{H}:=\sup\left\{p \in [0,1];
\mathbb{P}_{p}(|C_{x}|=\infty)=0\right\}$, we see that the
probability $\mathbb{P}_{p}(|C_{x}|=\infty)$ is equal to $0$, if
$p<p_{H}$ and strictly positive, if $p>p_{H}$. In the case
$p<p_{H}$ there is no infinite cluster almost surely, while in the
case $p>p_{H}$ there exists an infinite cluster almost surely. The
value of $p_{H}$ does not depend on the vertex $x$. It is often
called the \emph{percolation threshold}. The case $p<p_{H}$ is
called \emph{subcritical phase}, the case $p>p_{H}$ is called
\emph{supercritical phase}, while $p=p_{H}$ is called
\emph{critical phase}.

If $p_{H}<1$ it is obvious that $\mathbb{E}_{p}(|C_{x}|)=\infty$
for all $p>p_{H}$. The behavior of $\mathbb{E}_{p}(|C_{x}|)$ in
the subcritical phase is much more interesting. For edge percolation, the expectation
$\mathbb{E}_{p}(|C_{x}|)$ has value $1$ at $p=0$ and is infinite
at $p=1$. Similarly, for site percolation this expectation is zero
at $p=0$ and is infinite at $p=1$.
It can be shown that the expectation value is a non-decreasing function of
$p$. So, if we define $p_{T}:=\sup\left\{p \in [0,1],
\mathbb{E}_{p}(|C_{x}|)<\infty\right\}$, we see that
$\mathbb{E}_{p}(|C_{x}|)$ is finite, if $p<p_{T}$ and infinite if
$p>p_{T}$. Like the value of $p_{H}$, the value of $p_{T}$ is also
independent of the choice of vertex $x$.

The relation $p_{T} \leq p_{H}$ between the critical values is
obvious. Our goal is to prove equality of the two values.
Our main result is the following.

\begin{theorem}\label{thm: main}
For every quasi-transitive graph $G$ we have $p_{T}=p_{H}$.
\end{theorem}
As mentioned in the introduction, for general percolation models on the lattice the equality of
the two critical points was proven in \cite{AizenmanB-87}.
The method of proof was the use of differential inequalities for
certain (finite volume) order parameters. In \cite{AizenmanB-87} one can also find a discussion
of the use of such differential inequalities in other models of statistical physics.
Using a different method, sharpness of the phase transition for site percolation on quasi-transitive graphs with subexponential
growth was proven in \cite{Menshikov-86}, see Remark \ref{rem: Menshikov} below.

Similarly as in the lattice setting \cite{AizenmanB-87}, Theorem \ref{thm: main} holds also for
long range and oriented percolation on quasi-transitive graphs. To show this, one has only to modify certain steps
in the proof of the basic version of Theorem \ref{thm: main}. We present and explain these modification
in the last section of this paper.
\bigskip

It is well known that, in the subcritical phase on the lattice, the probabilities of
the events of the form $\left\{|C_{x}| \geq n\right\}$ decay exponentially in $n$. The same result
holds in the case of quasi-transitive graphs.

\begin{theorem}\label{thm: exponential decay}
Let $G$ be a quasi-transitive graph and let $p < p_{H}$. We can
find a constant $\alpha_{p} > 0$ such that for any positive
integer $n$ we have
$$
\mathbb{P}_{p}(|C_{x}| \geq n) \leq e^{-\alpha_{p}n}, \ \textrm{
for any vertex } x.
$$
\end{theorem}

In the lattice case, exponential decay was first proven for all
$p$ such that $\mathbb{E}_{p}(|C|) < \infty$. This result follows
from Theorem 5.1 in \cite{Kesten-82}. The same result was
proven for more general models on transitive graphs in Proposition
5.1 in \cite{AizenmanN-84}. Consequently, the exponential decay in the subcritical phase
is just a corollary of the equality of critical points
$p_{T}$ and $p_{H}$. The proof of Proposition 5.1 from \cite{AizenmanN-84}
extends directly from transitive graphs to quasi-transitive ones.
Thus Theorem \ref{thm: exponential decay} follows directly from Theorem \ref{thm: main}.

\begin{remark} \label{rem: Menshikov}
In \cite{Menshikov-86,MenshikovMS-86} Menshikov et al.~pursued a different route of argument.
They first show that for site percolation with $ p <p_H$ on quasi-transitive graphs of subexponential growth
the cluster radius distribution decays exponentially. More precisely, for every
$ p <p_H$ there exists a constant $\tilde\alpha_{p} > 0$ such that for all $x \in V$ and
all $n \in \NN$
\begin{equation}
\label{eq: radius distribution}
\mathbb{P}_{p}(C_{x} \cap S(x,n) \neq \emptyset ) \leq e^{-\tilde \alpha_{p}n}
\end{equation}
holds. By the subexponential growth condition on the graph, this
implies that the expected cluster size is finite. The key step in the proof of
\eqref{eq: radius distribution} is an estimate on the conditional expectation
\begin{equation*}
\mathbb{E}_{p}( \ | \delta \{C_{x} \cap S(x,n) \neq \emptyset\}| \ \mid \ C_{x} \cap S(x,n) \neq \emptyset ) ,
\end{equation*}
where $|\delta A|$ denotes the number of pivotal sites for the event $A$.
Note that the estimate \eqref{eq: radius distribution} on the cluster radius distribution is
weaker than the one in Theorem \ref{thm: exponential decay} on the cluster size distribution.
\end{remark}
\medskip

In the proof of Theorem \ref{thm: main} we will need the following
result.

\begin{proposition}\label{prop: mean cluster size}
For percolation on a quasi-transitive graph, we have for every vertex $x$:
\begin{equation}\label{eq: cluster size divergence}
\lim_{p \uparrow p_{T}}\mathbb{E}_{p}(|C_{x}|) = \infty.
\end{equation}
In particular, $\mathbb{E}_{p_{T}}(|C_{x}|) = \infty$ for any  $x\in V$.
\end{proposition}

In the lattice case the divergence of $\mathbb{E}_{p_{T}}(|C_{x}|)$ was proven in Corollary 5.1 in
\cite{Kesten-82}. The stronger statement \eqref{eq: cluster size divergence}
was then proven for more general percolation processes on transitive graphs in Lemma 3.1 in
\cite{AizenmanN-84}. The proof of this statement for quasi-transitive
graphs is essentially the same.

The versatility of the  differential inequalities method as presented in
\cite{AizenmanB-87} is illustrated by the fact that on the way to prove Theorem
\ref{thm: main} one obtains as an aside a bound on the critical exponent $\delta$,
cf.~\eqref{eqn: critical exponent} for a definition.

\begin{proposition}
The critical exponent satisfies $\delta \ge 2$.
\end{proposition}
This is a direct consequence of Lemma \ref{lemma: critical
asymptotics}.

\section{Basic facts}
\label{s-basics}

Now we shall present some basic definitions and results from
percolation theory. To be able to treat both the site and the bond model simultaneously, we shall denote, in the bond case, the edge set of a given graph $G$ by $S$. In the site case, $S$ will denote the vertex set of $G$.

\begin{definition}\label{def: events}
\begin{itemize}
\item[a)] We say that the event $A \in \mathcal{A}$ is
\emph{increasing}, if
$$
\omega_1 \in A, \ \ \omega_1 \leq \omega_{2} \  \Rightarrow \
\omega_{2} \in A.
$$
Here elements of $\Omega$ are ordered as functions from $S$ to
$\left\{0,1\right\}$. \item[b)] We say that a random variable $N$
is \emph{increasing}, if for any two configurations $\omega_{1}$
and $\omega_{2}$, such that $\omega_{1} \leq \omega_{2}$ we have
$N(\omega_{1}) \leq N(\omega_{2})$. \item[c)] We say that an event
$A$ depends only on finitely many states, if it is contained in
some finite dimensional cylinder set in $\mathcal{A}$. \item[d)]
For two increasing events $A_{1}$ and $A_{2}$, which depend only
on finitely many states, we define the event
\begin{multline*}
A_{1} \circ A_{2} :=\left\{\omega \in \Omega; \textrm{ there are
disjoint } S_{1}, S_{2} \subset \supp \omega, \textrm{ such that
for any } \omega_{1},\omega_{2} \in \Omega \right. \\ \left.
\omega_{i}|_{S_{i}}=1 \Rightarrow \omega_{i} \in A_{i}, \ i=1,2
\right\},
\end{multline*}
where $\supp \omega :=\left\{s \in S; \omega_{s}=1\right\}$.
\item[e)] For an increasing event $A$ and $\omega \in \Omega$ we
say that $s_{0} \in S$ is \emph{pivotal} for $A$ with respect to
$\omega$, if $\omega_{1} \in A$ and $\omega_{0} \notin A$, where
$\omega_{0}$ and $\omega_{1}$ have the same values as $\omega$ on
all elements of $S$ except on $s_{0}$ where $\omega_{i}$ has value $i$ ($i=0,1$). The
set $\left\{s_{0} \textrm{ is pivotal for the event } A\right\}$ is
obviously an event.
\end{itemize}
\end{definition}

\begin{tools}
\begin{itemize} \item[\upshape{a)}] For
any increasing event $A$ the function $p \mapsto
\mathbb{P}_{p}(A)$ is non-decreasing. \item[\upshape{b)}] For any
increasing random variable $N$, the function $p \mapsto
\mathbb{E}_{p}(N)$ is
non-decreasing. \item[\upshape{c)}] \textbf{Russo formula}\\
Suppose $A$ is an increasing event which depends only on states of
elements in some finite set $S'$, more precisely, on $\omega|_{S'}$, where
$S' \subset S$ is finite. Let $\textbf{p}=(\textbf{p}_{s})_{s \in
S'}$ be a given vector, such that $\textbf{p}_{s} \in [0,1]$, for
all $s \in S'$. Let $\mathbb{P}_{\textbf{p}}$ be the product probability
measure constructed in the same way as
the percolation measure before, by declaring an $s \in S'$ to
be open with probability $\textbf{p}_{s}$. Then the function
$\textbf{p} \mapsto \mathbb{P}_{\textbf{p}}$ has all first partial
derivatives, which satisfy
$$
\frac{d\mathbb{P}_{\textbf{p}}}{d\textbf{p}_{s}} (A) =
\mathbb{P}_{\textbf{p}}(s \textrm{ is pivotal for } A), \ \textrm{
for any } s \in S'.
$$
\item[\upshape{d)}] \textbf{FKG inequality}\\ For any increasing events
$A_{1}$ and $A_{2}$ we have $\mathbb{P}(A_{1} \cap A_{2}) \geq
\mathbb{P}(A_{1})P(A_{2})$. \item[\upshape{e)}] \textbf{BK inequality}\\ For
any increasing events $A_{1}$ and $A_{2}$, which depend only on
finitely many states, we have $\mathbb{P}(A_{1} \circ A_{2}) \leq
\mathbb{P}(A_{1})\mathbb{P}(A_{2})$.
\end{itemize}
\end{tools}

For the proofs of these Fundamental Tools and more background see Chapter 2 in
\cite{Grimmett-99}.

\begin{remark} \label{rem: vertex independence}
Having these results one can easily prove some claims from the previous section.
Since the event $\left\{|C_{x}|=\infty\right\}$ is increasing,
the function $p \mapsto \mathbb{P}_{p}(|C_{x}|=\infty)$ is non-decreasing.
Similarly, the random variable $|C_{x}|$ is increasing which implies that
the function $p \mapsto \mathbb{E}_{p}(|C_{x}|)$ is non-decreasing.
Using the FKG inequality one can easily show that the constants
$p_{T}$ and $p_{H}$ do not depend on the choice of the vertex $x$.
To see this, first notice
\begin{equation}\label{eq: vertex independance}
\left\{x \textrm{ is connected to } y\right\} \cap
\left\{|C_{y}|=\infty\right\} \subseteq
\left\{|C_{x}|=\infty\right\}.
\end{equation}
Now, because $G$ is connected, $\mathbb{P}_{p}(x \textrm{ is
connected to } y) > 0$ and thus the FKG inequality and (\ref{eq:
vertex independance}) imply 
\[
\mathbb{P}_{p}(|C_{y}|=\infty)>0
\Rightarrow \mathbb{P}_{p}(|C_{x}|=\infty)>0. 
\]
Because of the
symmetry of the role played by $x$ and $y$, $p_{H}$ does not
depend on $x$. To prove the claim for $p_{T}$, one decomposes
$\mathbb{E}_{p}(|C_{x}|)$ in the form $\displaystyle
\mathbb{E}_{p}(|C_{x}|) =
\sum_{n=1}^{\infty}\mathbb{P}_{p}(|C_{x}|\geq n)$ and uses a
relation similar to (\ref{eq: vertex independance}) for the
increasing event $\left\{|C_{x}|\geq n\right\}$.
\end{remark}

\section{The order parameter}
\label{s-order-parameter}
In this and the following section we shall work exclusively in the nearest neighbor Bernoulli bond percolation model. The site model will be discussed in section \ref{s-site-model}.

In the remainder of the paper it will be more convenient to work with
the parameter $\beta>0$ such that $p=1-e^{-\beta}$, instead with
$p$. Assuming $p_{H}<1$, we can define $\beta_{T}$ and $\beta_{H}$
by $p_{T}=1-e^{-\beta_{T}}$ and $p_{H}=1-e^{-\beta_{H}}$. We will
prove Theorem \ref{thm: main} in the context of $\beta_{T}$ and
$\beta_{H}$, but our proof works also in the case $p_{H}=1$ (this
case corresponds to $\beta_{H}=\infty$). Also we will abuse
notation by writing $\mathbb{P}_{\beta}$ for the probability
measure which corresponds to percolation with parameter
$p=1-e^{-\beta}$, and use a similar notation for the expectation.

From now on we will assume that we are given a fixed quasi-transitive
graph $G$. Subgraphs of $G$, which we will consider, will not be
required to be quasi-transitive. Moreover, we will assume that
some fundamental domain $\mathcal{F}$ is chosen and fixed.
For each positive integer $l$ we define a subgraph $\Lambda_{l}$
as follows. For the set of vertices of the graph $\Lambda_{l}$
take $\bigcup_{x \in \mathcal{F}}B(x,l)$ and connect two vertices,
if and only if they are connected in the graph $G$.
\medskip

To prove Theorem \ref{thm: main}, we will follow the arguments in
\cite{AizenmanB-87}. The idea of the proof there was to consider a so
called \emph{order parameter}, a function of two variables which contains
information about both $\mathbb{P}_{\beta}(|C_{x}|=\infty)$ and
$\mathbb{E}_{\beta}(|C_{x}|)$. For any vertex $y$ we define the
\emph{order parameter with respect to $y$} by
$$M_{y} \colon \left]0,\infty\right[^{2} \to
[0,1], \ \ \ M_{y}(\beta,h):=1-\sum_{n \in
\mathbb{N}}\mathbb{P}_{\beta}(|C_{y}|= n)e^{-nh}.
$$
The \emph{order parameter} $M$ is defined as
$$M\colon ]0,\infty[^{2}
\to \mathbb{R}, \ \ \ M(\beta,h):=\sum_{x \in
\mathcal{F}}M_{x}(\beta,h).
$$
When a finite subgraph $G'$ is given, we can define an analog
function with respect to $G'$. Namely, for $y$ in $G'$ we define
$\displaystyle M_{y}^{G'}(\beta,h):=1-\sum_{n \in
\mathbb{N}}\mathbb{P}(|C_{y}^{G'}| = n)e^{-nh}$. Particularly
interesting for our purposes will be the \emph{finite volume order
parameter}, defined as $\displaystyle
M^{\Lambda_{l}}(\beta,h):=\sum_{x \in
\mathcal{F}}M_{x}^{\Lambda_{l}}(\beta,h)$.

In the following Lemmata and Propositions we establish certain
basic properties of the order parameter.

\begin{lemma}\label{lemma: OP definition}
Let $G'$ be an arbitrary subgraph of $G$ and $y$ an arbitrary
vertex in $G'$. The following formula holds
\begin{equation}\label{eq: OP definition}
M_{y}^{G'}(\beta,h)=\sum_{n \in \mathbb{N}}\mathbb{P}(|C_{y}^{G'}|
\geq n)(e^{-(n-1)h}-e^{-nh}).
\end{equation}
In particular, {\upshape(\ref{eq: OP definition})} holds in the
cases $G'=G$ and $G'=\Lambda_{l}$.
\end{lemma}
\begin{proof}
The proof is straightforward, using $\displaystyle
\mathbb{P}(|C_{y}^{G'}| \geq
n)=1-\sum_{k=1}^{n-1}\mathbb{P}(|C_{y}^{G'}|=k)$.
\end{proof}

\begin{proposition}\label{prop: OP properties}
The order parameter $M$ has the following properties.
\begin{itemize}
\item[\upshape{a)}] $M$ is a non-decreasing function in both
variables. \item[\upshape{b)}] $M$ has a continuous partial
derivative in $h$, and we have the formula
\begin{equation}\label{eq: OP properties}
\frac{\partial M}{\partial h}(\beta,h) = \sum_{x \in
\mathcal{F}}\sum_{n \in
\mathbb{N}}n\mathbb{P}_{\beta}(|C_{x}|=n)e^{-nh}.
\end{equation}
\end{itemize}
The analogous claims hold for the finite volume order parameter.
\end{proposition}

\begin{proof}
\begin{itemize}
\item[a)] Since the event $\left\{|C_{y}| \geq n\right\}$ is an
increasing event, the probability $\mathbb{P}_{\beta}(|C_{y}| \geq
n)$ is a non-decreasing function of $\beta$. From Lemma
\ref{lemma: OP definition} it is clear that $M$ is non-decreasing
in $\beta$. On the other hand, from the definition it is clear
that $M$ is even strictly increasing in $h$. \item[b)] To prove
this claim we just have to show that the series of formal partial
derivatives $\displaystyle \sum_{n \in
\mathbb{N}}n\mathbb{P}_{\beta}(|C_{y}|=n)e^{-nh}$ converges
locally uniformly. But this is clear since
$$
\sum_{n \in \mathbb{N}}n\mathbb{P}_{\beta}(|C_{y}|=n)e^{-nh} \leq
\sum_{n \in \mathbb{N}}ne^{-nh}
$$
and the latter series converges absolutely and locally uniformly.
\end{itemize}
\end{proof}

The functions $M$ and $\frac{\partial M}{\partial h}$ are positive
on $]0,\infty[^{2}$, $h \mapsto M(\beta,h)$ is non-decreasing and
$h \mapsto \frac{\partial M}{\partial h}(\beta,h)$ is
non-increasing. The last claim is clear from the formula (\ref{eq:
OP properties}). Thus the limits $\lim_{h\downarrow 0}M(\beta,h)$ and $\lim_{h \downarrow 0}\frac{\partial
M}{\partial h}(\beta,h)$ are well-defined with values in $[0,\infty[$, respectively
$[0,\infty]$. The next proposition gives the
probabilistic interpretation of these limits.

\begin{proposition}\label{prop: probability limits}
For every $\beta \in ]0,\infty[$ we have the following
\begin{equation}\label{eq: probability limits}
\lim_{h \downarrow 0}M(\beta,h) = \sum_{x \in
\mathcal{F}}\mathbb{P}_{\beta}(|C_{x}|=\infty), \ \ \ \ \lim_{h
\downarrow 0}\frac{\partial M}{\partial h}(\beta,h) = \sum_{x \in
\mathcal{F}}\mathbb{E_{\beta}}(|C_{x}|;|C_{x}|<\infty).
\end{equation}
\end{proposition}

\begin{proof}
Since $\lim_{h \downarrow 0}e^{-nh}=1$, for every $n \in
\mathbb{N}$, using the Lebesgue monotone convergence theorem, we
get $\displaystyle \lim_{h \downarrow 0}M_{x}(\beta,h) = 1-\sum_{n
\in \mathbb{N}}\mathbb{P}_{\beta}(|C_{x}|=n)=
\mathbb{P}_{\beta}(|C_{x}|=\infty)$. Now the first
equality in (\ref{eq: probability limits}) follows. The second one can be
proved in the same manner, using formula (\ref{eq: OP
properties}).
\end{proof}

Thus we obtain --- as indicated earlier --- the functions
$ \beta \mapsto \sum_{x \in \mathcal{F}}\mathbb{P}_{\beta}(|C_{x}|=\infty)$
and $ \beta \mapsto \sum_{x \in \mathcal{F}}\mathbb{E_{\beta}}(|C_{x}|;|C_{x}|<\infty)$,
which we wanted to understand in the first place, as marginals
of $M$ and $ \frac{\partial M}{\partial h}$.
Now we give two lemmata, which will be used repeatedly in the
proof of the key inequalities presented in Propositions
\ref{prop: inequality 1} and \ref{prop: inequality 2}.

\begin{lemma}\label{lemma: cluster size inequalities}
\begin{itemize}
\item[\upshape{a)}] Let $G_{1}$ be any subgraph of $G$, and
$G_{2}$ any finite subgraph of $G_{1}$, containing some vertex
$x$. For any nonnegative integer $n$ we have
\begin{equation}\label{eq: cluster size inequalities 0}
\mathbb{P}(|C_{x}^{G_{2}}| \geq n) \leq \mathbb{P}(|C_{x}^{G_{1}}|
\geq n).
\end{equation}
Moreover, $M_{x}^{G_{2}}(\beta,h) \leq M_{x}^{G_{1}}(\beta,h)$,
for all positive $\beta$ and $h$.
\item[\upshape{b)}] Let $y$ be a vertex of $\Lambda_{l}$ and
$x$ be the unique element of $\mathcal{F}$ in the same
orbit as $y$. For any nonnegative integer $n$ we have
\begin{equation}\label{eq: nejednakosti 1}
\mathbb{P}(|C_{y}^{\Lambda_{l}}| \geq n)
\leq \mathbb{P}(|C_{y}| \geq n) = \mathbb{P}(|C_{x}|\geq n).
\end{equation}
\item[\upshape{c)}] For any vertex $x \in \mathcal{F}$ and $n \leq
l$ we have
\begin{equation}\label{eq: nejednakosti 2}
\mathbb{P}(|C_{x}^{\Lambda_{l}}| \geq n) = \mathbb{P}(|C_{x}| \geq
n).
\end{equation}
\end{itemize}
\end{lemma}

\begin{proof}
\begin{itemize}
\item[a)] Let $A$ be an arbitrary connected subgraph of $G_{2}$,
containing the vertex $x$. The identification from Remark \ref{rem: 1}
implies that the probabilities of the events
$\left\{C_{x}^{G_{2}}=A\right\}$ and $\left\{A \textrm{ is the
component of }C_{x} \cap G_{2} \textrm{ containing } x\right\}$
are equal. Similarly the probability of the event $\left\{A
\textrm{ is the component of } C_{x}^{G_{1}} \cap G_{2} \textrm{
containing } x \right\}$ is equal to probability of $\left\{A
\textrm{ is the component of } C_{x} \cap G_{2} \textrm{
containing } x\right\}$. So we can write
$$\mathbb{P}(C_{x}^{G_{2}}=A) =
\mathbb{P}(A \textrm{ is the component of } C_{x}^{G_{1}} \cap
G_{2} \textrm{ containing } x).
$$
Since the events on the right side are disjoint for different
$A$'s, we can write
\begin{align} \nonumber
\mathbb{P}(|C_{x}^{G_{2}}| \geq n) &
= \sum_{A; |A|\geq n}\mathbb{P}(C_{x}^{G_{2}} = A) 
\\ \nonumber
&= \sum_{A;|A|\geq n} \mathbb{P}(A\textrm{ is the component of } C_{x}^{G_{1}} \cap G_{2} 
\text{ containing } x)
\\ \label{eq: cluster size inequalities 1}
&\le  \mathbb{P}(|C_{x}^{G_{1}} \cap G_{2}| \geq n)
\leq \mathbb{P}(|C_{x}^{G_{1}}| \geq n).
\end{align}
The last inequality follows from $\left\{|C_{x}^{G_{1}} \cap
G_{2}| \geq n\right\} \subset \left\{|C_{x}^{G_{1}}| \geq
n\right\}$. The sums in (\ref{eq: cluster size inequalities 1})
are taken over all connected subgraphs of $G_{2}$, which contain
$x$ and which are of size greater or equal than $n$. Since $G_{2}$
is a finite graph, these sums are finite. Now the claim
$M_{x}^{G_{2}}(\beta,h) \leq M_{x}^{G_{1}}(\beta,h)$ follows
directly from Lemma \ref{lemma: OP definition}. \item[b)] The
first inequality follows directly from part a), if we take $G_2=
\Lambda_l$ and $G_1= G$. The second (in)equality follows from the
fact that there is an automorphism $\tau$ such that $x=\tau y$
and that the probability measure is invariant under $\tau$.
 \item[c)] For any $k<l$ and any connected subgraph $A$ of
$\Lambda_{l}$ of size $k$, which contains $x$, the edge set and
the edge boundary of $A$ are contained in $\Lambda_{l}$. So it is
clear that for any such $A$ we have
$\mathbb{P}(C_{x}^{\Lambda_{l}}=A) = \mathbb{P}(C_{x}=A)$. Taking
the sum over all possible $A$'s, when $k$ is fixed, we obtain
$\mathbb{P}(|C_{x}^{\Lambda_{l}}| = k) = \mathbb{P}(|C_{x}| = k)$.
Taking the sum over $k<n$ and subsequently complements yields the
result.
\end{itemize}
\end{proof}

\begin{lemma}\label{lemma: basic modification}
Let $y$ be a vertex of $\Lambda_{l}$ and $x$ the unique element of
$\mathcal{F}$ in the same orbit as $y$. Then, for all $(\beta,h) \in
]0,\infty[^{2}$, the following inequality holds
\begin{equation}\label{eq: basic modification}
M_{y}^{\Lambda_{l}}(\beta,h) \leq M_{x}^{\Lambda_{l}}(\beta,h) +
e^{-lh}.
\end{equation}
\end{lemma}

\begin{proof}
Using Lemma \ref{lemma: OP definition} we can write
\begin{equation}\label{eq: basic modification proof 1}
M_{y}^{\Lambda_{l}}(\beta,h) =
\sum_{n=1}^{l}\mathbb{P}(|C_{y}^{\Lambda_{l}}|\geq
n)(e^{-(n-1)h}-e^{-nh}) + \sum_{l+1 \leq n <
\infty}\mathbb{P}(|C_{y}^{\Lambda_{l}}|\geq
n)(e^{-(n-1)h}-e^{-nh}).
\end{equation}
Using parts b) and c) of Lemma \ref{lemma: cluster size
inequalities} we can bound the first summand
\begin{multline}\label{eq: basic modification proof 2}
\sum_{n=1}^{l}\mathbb{P}(|C_{y}^{\Lambda_{l}}|\geq
n)(e^{-(n-1)h}-e^{-nh}) \leq \sum_{n=1}^{l}\mathbb{P}(|C_{x}|\geq
n)(e^{-(n-1)h}-e^{-nh})
\\ = \sum_{n=1}^{l}\mathbb{P}(|C_{x}^{\Lambda_{l}}|\geq
n)(e^{-(n-1)h}-e^{-nh}) \leq M_{x}^{\Lambda_{l}}(\beta,h).
\end{multline}
The second summand can be easily bounded
\begin{equation}\label{eq: basic modification proof 3}
\sum_{l+1 \leq n < \infty}\mathbb{P}(|C_{y}^{\Lambda_{l}}|\geq
n)(e^{-(n-1)h}-e^{-nh}) \leq \sum_{l+1 \leq n <
\infty}(e^{-(n-1)h}-e^{-nh}) = e^{-lh}.
\end{equation}
Inserting (\ref{eq: basic modification proof 2}) and (\ref{eq:
basic modification proof 3}) into (\ref{eq: basic modification
proof 1}) we get the result.
\end{proof}

\begin{remark}
In \cite{AizenmanB-87}, where percolation on the lattice $\ZZ^d$
was analyzed, the finite graphs $\Lambda_l$ where chosen to be tori,
i.e. cubes with periodic boundary conditions. This has the advantage
that percolation on the finite graphs is still homogeneous under translations.
In this situation, \eqref{eq: basic modification} simplifies to
$M_{y}^{\Lambda_{l}}(\beta,h) = M_{x}^{\Lambda_{l}}(\beta,h)$.
In particular, there is no finite volume correction term $e^{-lh}$.
\end{remark}

Notice that the subgraphs $\Lambda_{l}$ exhaust the whole graph $G$ as
$l$ goes to $\infty$. Therefore, in the macroscopic limit $l \to \infty$, we
can expect the finite volume order parameters to behave like the
order parameter. The proof of this claim is the content of the
next proposition.

\begin{proposition}\label{lemma: convergence}
The finite volume order parameter $M^{\Lambda_{l}}$ and its
partial derivative $\frac{\partial M^{\Lambda_{l}}}{\partial h}$
converge pointwise to $M$ and $\frac{\partial M}{\partial h}$,
respectively. In other words, 	for all $(\beta,h)\in]0,\infty[^{2}$ we have
$$
\lim_{l \to \infty} M^{\Lambda_{l}}(\beta,h) = M(\beta,h), \ \ \
\lim_{l \to \infty} \frac{\partial M^{\Lambda_{l}}}{\partial
h}(\beta,h) = \frac{\partial M}{\partial h}(\beta,h).
$$
\end{proposition}

\begin{proof}
Fix some $x \in \mathcal{F}$. Lemma \ref{lemma: cluster size
inequalities} c) implies $\lim_{l \to \infty}
\mathbb{P}(|C_{x}^{\Lambda_{l}}| \geq n) = \mathbb{P}(|C_{x}| \geq
n)$, for any positive integer $n$. Using Lemma \ref{lemma: OP
definition} and the Lebesgue dominated convergence theorem we get
$$
\lim_{l \to \infty} M_{x}^{\Lambda_{l}}(\beta,h) = M_{x}(\beta,h),
\ \ \textrm{ for all } (\beta,h)\in ]0,\infty[^{2},
$$
for any $x \in \mathcal{F}$. Taking the sum over $x \in
\mathcal{F}$ we get the desired result for the order parameter.
The claim for the partial derivative is obtained in the same way
using the formula
$$\frac{\partial M^{\Lambda_{l}}}{\partial h} =
\sum_{x \in \mathcal{F}}\sum_{n \in
\mathbb{N}}n\mathbb{P}_{\beta}(|C_{x}^{\Lambda_{l}}|=n)e^{-nh},$$
which is proved in the same way as the formula for $\frac{\partial
M}{\partial h}$ in Proposition \ref{prop: OP properties} b).
\end{proof}

Another way of looking at the order parameter is through the idea
of "colored sites", which was used in the paper \cite{AizenmanB-87}.
Fix a positive real $h>0$. For every vertex $y$ say that it
is \emph{blue} with  probability $1-e^{-h}$ independently of
all other vertices. The corresponding probability space is defined
similarly as the site percolation probability space before. For each vertex
$y$ define the probability space
$(\Omega_{y},P(\Omega_{y}),\mathbb{P}_{y})$, where
$\Omega_{y}:=\left\{0,1\right\}$, $P(\Omega_{y})$ is the power set
of $\Omega_{y}$ and $\mathbb{P}_{y}(1)=1-e^{-h}$,
$\mathbb{P}_{y}(0)=e^{-h}$. The probability space
$(\Omega',\mathcal{A}',\mathbb{P}_{h}')$ is defined as the product
of these probability spaces. From now on, we shall actually work
on the probability space $(\Omega,\mathcal{A},\mathbb{P}_{\beta})
\times (\Omega',\mathcal{A}',\mathbb{P}_{h}')$. The probability
measure will be denoted by $\mathbb{P}_{\beta,h}$, but again we
will often omit the subscript. The random set of blue sites will
be denoted by $B$. Analog functions can be defined on subgraphs of
$G$, and in this case we use the same notation as before.

The event that some vertex $y$ is connected to some blue site with
an open path will be denoted by $\left\{y \leftrightarrow
B\right\}$ or by $\left\{C_{y} \cap B \neq \emptyset\right\}$,
while $\left\{y \nleftrightarrow B\right\}$ and $\left\{C_{y} \cap
B = \emptyset\right\}$ will stand for the complement of this
event. For connectedness with open edges which have end-vertices
in some given set $A$ we will use $\left\{y \leftrightarrow_{A}
B\right\}$, while for connectedness with open edges which are also
edges of some subgraph $G'=(V',E')$ we will write $\left\{y
\leftrightarrow_{G'} B\right\}$. If $V'$ is a given finite set of
vertices and $A : \Omega \to P(V')$ a function from $\Omega$ to
the power set of $V'$ then $\left\{y \leftrightarrow _{A}
B\right\}$ will denote the set of $\omega$'s for which there is an
open path which connects $y$ with an element in $B$ and contains
only vertices in $A(\omega)$. If $A$ is such that for any subset
$V'' \subset V'$, the set $\left\{A(\omega) = V''\right\}$ is an
event, the set $\left\{y \leftrightarrow _{A} B\right\}$ is also
an event. Similar notation will be used for random subgraphs.

The blue sites are, in some sense, identified with "infinity". For
example, for some fixed vertex $x$, the event $\left\{x
\leftrightarrow B \right\}$ is the generalization of
$\left\{|C_{x}|=\infty\right\}$, because the open path from $x$
which reaches some vertex in $B$ is considered to have escaped to
infinity. Intuitively, if the parameter $h$ decreases to $0$, the
density of blue sites decreases to $0$, and they "move further and
further away" from $x$. So, their effect on the whole picture gets
less relevant and in the limit $h \downarrow 0$ we should expect
to return to our original percolation setting. Namely, the
probability of the event $\left\{x \leftrightarrow B \right\}$
should converge to the probability of the event
$\left\{|C_{x}|=\infty\right\}$. This is actually a direct consequence of
Proposition \ref{prop: probability limits} in view of Proposition
\ref{prop: blue sites}. The next proposition shows the relationship between the
order parameter and blue sites.

\begin{proposition}\label{prop: blue sites}
Let $G' \subset G$ be a subgraph of a quasi-transitive graph $G$
and $y$ some vertex in $G'$. Using the above notation we have
\begin{equation}\label{eq: blue sites}
M_{y}(\beta,h) = \mathbb{P}_{\beta,h}(y \leftrightarrow B), \ \ \
M_{y}^{G'}(\beta,h) = \mathbb{P}_{\beta,h}(y \leftrightarrow _{G'}
B).
\end{equation}
\end{proposition}

\begin{proof}
We prove the first equality in (\ref{eq: blue sites}), while the
second one can be proven in the same way. It is enough to show
\begin{equation}\label{eq: blue sites proof 0}
\mathbb{P}(y \nleftrightarrow B) = \sum_{n \in
\mathbb{N}}\mathbb{P}(|C_{y}|=n)e^{-nh}.
\end{equation}
For any positive integer $n$ we have $$\mathbb{P}(|C_{y}|=n,y
\nleftrightarrow B) = \sum_{A;|A|=n}\mathbb{P}(C_{y}=A,A \cap B
=\emptyset) = \sum_{A;|A|=n}\mathbb{P}(C_{y}=A)\mathbb{P}(A \cap B
= \emptyset),$$ where the last equality is obtained using the
independence of bond and site variables and the sums are taken
over all connected subgraphs with $n$ vertices containing $y$.
Obviously $\mathbb{P}(A \cap B = \emptyset)=e^{-h|A|}$ and so one
obtains,
\begin{equation}\label{eq: blue sites proof 1}
\mathbb{P}(|C_{y}|=n,y \nleftrightarrow B) =
\mathbb{P}(|C_{y}|=n)e^{-nh}.
\end{equation}
Now, we are left to estimate $\mathbb{P}(|C_{y}|=\infty ,y
\nleftrightarrow B)$. Define the random variable
$k_{n}:=\min\left\{m;|C_{y}\cap B(y,m)|\geq n\right\}$ which
obviously has only finite values on the event
$\left\{|C_{y}|=\infty\right\}$. Next we can write
\begin{multline}\label{eq: blue site proof 1-1}
\mathbb{P}(k_{n}=m,y\nleftrightarrow B) \leq
\sum_{A}\mathbb{P}(C_{y}\cap B(y,m)=A, A\cap B = \emptyset) \\
\leq \sum_{A}\mathbb{P}(C_{y}\cap B(y,m)=A)e^{-h|A|} \leq
e^{-hn}\mathbb{P}(k_{n}=m),
\end{multline}
where the sum is taken over all possible realization of $C_{y}\cap
B(y,m)=A$ such that the condition $k_{n}=m$ is fulfilled. Since
$$
\mathbb{P}(|C_{y}|=\infty ,y \nleftrightarrow B) \leq
\mathbb{P}(k_{n}<\infty ,y \nleftrightarrow B)=
\sum_{m=0}^{\infty}\mathbb{P}(k_{n}=m,y \nleftrightarrow B),
$$
for all positive integers $n$, using (\ref{eq: blue site proof
1-1}) we obtain $\mathbb{P}(|C_{y}|=\infty ,y \nleftrightarrow B)
\leq e^{-hn}$ for all positive integers $n$ and thus
\begin{equation}\label{eq: blue sites proof 1-2}
\mathbb{P}(|C_{y}|=\infty ,y \nleftrightarrow B) =0.
\end{equation}
From (\ref{eq: blue sites proof 1}) and (\ref{eq: blue sites proof
1-2}) we get (\ref{eq: blue sites proof 0}), and hence the proof is
completed.
\end{proof}

\section{Differential inequalities for the order parameter}
\label{s-diffinequalities}

In the following we will prove two differential inequalities
involving the order parameter $M^{\Lambda_{l}}$. These
inequalities differ from the inequalities (3.1) and (3.2) in
\cite{AizenmanB-87}, by the additional term $e^{-hl}$.
This finite volume correction appears in our situation,
since for general automorphism groups we cannot use ``periodic''
boundary conditions and thus percolation on the finite graph $\Lambda_l$
is no longer homogeneous. For residually finite automorphism groups it may
be possible to use periodic boundary conditions and thus eliminate the
correction term $e^{-hl}$.

The abovementioned differential inequalities will be crucial for the proof of Theorem \ref{thm: main},
because they contain essential information about the behavior of the order
parameter $M(\beta,h)$ when $h$ approaches $0$. In the proof of
Theorem \ref{thm: main} we will forget about the percolation and probability
setting, and work with these inequalities instead,
using analytic methods. As for the proof of these inequalities, we
will use the notion of blue sites extensively.

\begin{proposition}\label{prop: inequality 1}
There exists a constant $K>0$ such that
\begin{equation}\label{eq: inequalities 1}
\frac{\partial M^{\Lambda_{l}}}{\partial \beta} \leq K
(M^{\Lambda_{l}}+e^{-lh})\frac{\partial M^{\Lambda_{l}}}{\partial
h}, \textrm{ for all positive } \beta \textrm{ and } h.
\end{equation}
\end{proposition}

\begin{proof}
The event $\left\{x \leftrightarrow
_{\Lambda_{l}}B\right\}$ is increasing and depends on the states
of only finitely many edges. Using the Russo formula for this
event, we obtain
$$
\frac{\partial M_{x}^{\Lambda_{l}}}{\partial \beta} = \sum_{[y,z]
\in \Lambda_{l}}e^{-\beta}\mathbb{P}([y,z] \textrm{ is pivotal for
the event } \left\{ x \leftrightarrow _{\Lambda_{l}}B \right\}),
$$
where the sum is taken over all edges $[y,z]$ in $\Lambda_{l}$.
Since $\mathbb{P}([y,z] \textrm{ is closed})=e^{-\beta}$, and the
events $\left\{[y,z] \textrm{ is closed}\right\}$ and
$\left\{[y,z] \textrm{ is pivotal for the event } \left\{ x
\leftrightarrow _{\Lambda_{l}}B \right\}\right\}$ are independent,
we get
$$
\frac{\partial M_{x}^{\Lambda_{l}}}{\partial \beta} = \sum_{[y,z]
\in \Lambda_{l}} \mathbb{P}([y,z] \textrm{ is closed},[y,z]
\textrm{ is pivotal for the event } \left\{ x \leftrightarrow
_{\Lambda_{l}}B \right\}).
$$
One should notice that
\begin{multline*}
\left\{[y,z] \textrm{ is closed},[y,z] \textrm{ is pivotal for the
event } \left\{ x \leftrightarrow _{\Lambda_{l}}B
\right\}\right\}\\ = \left\{ x \nleftrightarrow _{\Lambda_{l}}B,
[y,z] \textrm{ is pivotal for the event } \left\{ x
\leftrightarrow _{\Lambda_{l}}B
\right\}\right\} \\
= \left\{ C_{x}^{\Lambda_{l}} \cap B = \emptyset, y \in
C_{x}^{\Lambda_{l}}, z \leftrightarrow _{\Lambda_{l}\backslash
C_{x}^{\Lambda_{l}}}B\right\} \bigcup \left\{C_{x}^{\Lambda_{l}}
\cap B = \emptyset, z \in C_{x}^{\Lambda_{l}}, y \leftrightarrow
_{\Lambda_{l}\backslash C_{x}^{\Lambda_{l}}}B\right\}.
\end{multline*}
Here $\Lambda_{l}\backslash C_{x}^{\Lambda_{l}}$ stands for the
graph obtained by deleting the vertices in $C_{x}^{\Lambda_{l}}$
and all incident edges of the graph $\Lambda_{l}$. Similar
notation will be used often in the rest of the paper. Now we pass
from a sum over undirected edges to a sum over directed ones and
write
\begin{align}\label{eq: inequalities proof 1}
\frac{\partial M_{x}^{\Lambda_{l}}}{\partial \beta} &= \sum_{(y,z)
\in \Lambda_{l}^{2} \atop y \sim z} \mathbb{P}(C_{x}^{\Lambda_{l}}
\cap B = \emptyset, y \in C_{x}^{\Lambda_{l}}, z \leftrightarrow
_{\Lambda_{l}\backslash C_{x}^{\Lambda_{l}}}B)\nonumber \\ & =
\sum_{(y,z) \in \Lambda_{l}^{2} \atop y \sim z} \sum_{A; y \in
A}\mathbb{P}(C_{x}^{\Lambda_{l}}=A, A\cap B = \emptyset, z
\leftrightarrow _{\Lambda_{l}\backslash A}B),
\end{align}
where the last sum is taken over all connected subgraphs $A$ of
$\Lambda_{l}$ containing $x$ and $y$. The event
$\left\{C_{x}^{\Lambda_{l}}=A\right\}$ depends only on the states
of edges which have at least one end-vertex in $A$. The event
$\left\{z \leftrightarrow _{\Lambda_{l}\backslash A}B\right\}$
depends only on the states of edges which do not have end-vertices
in $A$ and on the states of vertices outside $A$. Finally the
event $\left\{A \cap B = \emptyset\right\}$ depends only on the
states of vertices in $A$. Hence these events are
independent. Using this independence and Proposition \ref{prop:
blue sites}, equation (\ref{eq: inequalities proof 1}) can be rewritten as
\begin{equation}\label{eq: inequalities 1-2}
\frac{\partial M_{x}^{\Lambda_{l}}}{\partial \beta} = \sum_{(y,z)
\in \Lambda_{l}^{2} \atop y \sim z} \sum_{A; y \in
A}\mathbb{P}(C_{x}^{\Lambda_{l}}=A) \, M_{z}^{\Lambda_{l}\backslash A}
\, e^{-h|A|}.
\end{equation}
Lemmata \ref{lemma: cluster size inequalities} a) and \ref{lemma:
basic modification} imply $M_{z}^{\Lambda_{l}\backslash A} \leq
M_{z}^{\Lambda_{l}} \leq M^{\Lambda_{l}} + e^{-lh}$. Inserting
this into (\ref{eq: inequalities 1-2}) we obtain
\begin{align}\label{eq: inequalities proof 2}
\frac{\partial M_{x}^{\Lambda_{l}}}{\partial \beta} & \leq
K(M^{\Lambda_{l}}+e^{-lh}) \sum_{A} \sum_{y;y \in
A}\mathbb{P}(C_{x}^{\Lambda_{l}}=A)e^{-h|A|}\nonumber \\ & =
K(M^{\Lambda_{l}}+e^{-lh}) \sum_{A}|A| \,
\mathbb{P}(C_{x}^{\Lambda_{l}}=A) \, e^{-h|A|},
\end{align}
where the sum is taken over all possible graphs $A$ for
$C_{x}^{\Lambda_{l}}$ and $K$ is the maximal vertex degree in the
graph $G$. Now grouping together all $A$'s for which $|A|=n$, we
get
\begin{equation}\label{eq: inequalities proof 2-1}
\sum_{A}|A|\mathbb{P}(C_{x}^{\Lambda_{l}}=A)e^{-h|A|} = \sum_{n
\in \mathbb{N}}n\mathbb{P}(|C_{x}^{\Lambda_{l}}|=n)e^{-nh} =
\frac{\partial M_{x}^{\Lambda_{l}}}{\partial h}.
\end{equation}
Inserting (\ref{eq: inequalities proof 2-1}) into (\ref{eq:
inequalities proof 2}) and taking the sum over $x \in \mathcal{F}$ we
have proven the proposition.
\end{proof}

\begin{proposition}\label{prop: inequality 2}
The finite volume order parameter satisfies the following inequality
\begin{equation}\label{eq: inequalities 2}
M^{\Lambda_{l}} \leq h \frac{\partial M^{\Lambda_{l}}}{\partial h}
+ \big(M^{\Lambda_{l}}\big)^{2} + \beta
(M^{\Lambda_{l}}+e^{-lh})\frac{\partial M^{\Lambda_{l}}}{\partial
\beta}, \textrm{ for all positive } \beta \textrm{ and } h.
\end{equation}
\end{proposition}

\begin{proof}
In the proof of this inequality and especially in summations, $A$
will denote vertex sets.\\
To prove \eqref{eq: inequalities 2}, we have to change both our graph and
probability space. Let $n$ be an arbitrary, but fixed positive
integer. For every pair of adjacent vertices $y \sim z$ in
$\Lambda_{l}$ we replace the edge $[y,z]$ with $n$ edges which
will be denoted by $[y,z]_{1},[y,z]_{2},\dots [y,z]_{n}$. In this
way we obtain a new graph $G_{n}'$. We shall consider bond
percolation on the graph $G_{n}'$, and so we define the canonical
percolation product probability space in the usual way. A cluster
containing some vertex $x$ will be denoted by $C_{x}^{G_{n}'}$ and
its vertex set by $V_{x}^{G_{n}'}$. Notice that, for percolation
on $\Lambda_{l}$ with percolation parameter $1-e^{-\beta}$ and
for percolation on $G_{n}'$ with percolation parameter
$1-e^{-\beta /n}$, the probabilities that two adjacent vertices
are directly connected in the percolation graph are the same. This
implies the fact that the probability of the event
$\left\{V_{x}^{G_{n}'} = A\right\}$ in the new probability space
is equal to the probability of the event
$\left\{V_{x}^{\Lambda_{l}}=A\right\}$ in the old probability
space, for all possible sets of vertices $A$, if the parameter is
changed from $\beta$ to $\beta / n$. Here $V_{x}^{\Lambda_{l}}$
stands for the vertex set of $C_{x}^{\Lambda_{l}}$. Next we define
a graph $G_{n}$, which  contains $G_{n}'$ as a subgraph, by
adding to $G_{n}'$ a new vertex $b$, which is connected to
each of the vertices of $G_{n}'$ with exactly $n$ edges. The role
of the blue sites will be played by the edges incident to $b$. So,
in addition to the percolation on $G_{n}'$, we have the following rule:
We fix $h'
> 0$ and for every edge incident to $b$ we say that it is open with
probability $1-e^{-h'}$, independently of the states of all other
edges in the graph $G_{n}$. Notice that the events that some
vertex is blue, in the old probability space, and that some vertex
is directly connected to $b$, in the new probability space, have
the same probabilities, if $h'=h/n$ holds for the respective
parameters. This implies that the probabilities of the events
$\left\{x \leftrightarrow_{\Lambda_{l}} B\right\}$ and $\left\{x
\leftrightarrow_{G_{n}} b\right\}$ are the same, if both parameters
$\beta$ and $h$ are divided by $n$.
We will abbreviate the notation for the event $\left\{x\leftrightarrow_{G_{n}} b\right\}$
by writing simply by $\left\{x \leftrightarrow b\right\}$. So from now on we will
assume that $\beta$ and $h$ are fixed and we shall work in the new
probability space with parameters $\beta /n$ and $h/n$, and we
will keep in mind that the probabilities of the events mentioned
above remain unchanged. The probability measure will still be
denoted with $\mathbb{P}$.

We define the sets $F_{i}$, $i=1,2,3$, as follows:
\begin{align*}
&F_{1}:=\left\{\textrm{There is a unique open edge of } G_{n}
\textrm{ which connects some vertex of } C_{x}^{G_{n}'} \textrm{
with } b\right\},
\\ &
F_{2}:=\left\{x \leftrightarrow b\right\} \circ \left\{x
\leftrightarrow b\right\} = \left\{\textrm{There are two edge
disjoint paths from } x \textrm{ to } b\right\},
\\ &
F_{3}:=\bigcup_{(y,z) \in \Lambda_{l}^{2} \atop y \sim z}
\bigcup_{i=1}^{n}\left\{[y,z]_{i} \textrm{ is open and pivotal for
} \left\{x \leftrightarrow b\right\},\right. \\& \ \ \ \ \ \ \ \ \ \ \ \ \ \left.\left\{z
\leftrightarrow_{G_{n}\backslash [y,z]_{i}}b\right\} \circ
\left\{z \leftrightarrow_{G_{n}\backslash
[y,z]_{i}}b\right\}\right\},
\end{align*}
where $G_{n}\backslash [y,z]_{i}$ denotes the graph obtained from
$G_{n}$ by deleting the edge $[y,z]_{i}$. It is easy to see that
these sets are events. Lemma 3.5 from \cite{AizenmanB-87} implies that
$\left\{x \leftrightarrow b\right\}$ is a disjoint union of the
$F_{i}$'s and so $\displaystyle M_{x}^{\Lambda_{l}} = \mathbb{P}(x
\leftrightarrow b) = \mathbb{P}(F_{1}) + \mathbb{P}(F_{2}) +
\mathbb{P}(F_{3})$.

The probability $\mathbb{P}(F_{1})$ can be calculated as follows
\begin{align}\label{eq: inequalities proof 5}
\mathbb{P}(F_{1}) &= \sum_{A} \mathbb{P}(V_{x}^{G_{n}'}=A, A
\textrm{ is directly connected to } b \textrm{ with a unique open edge}) \nonumber \\
& = \sum_{A} \mathbb{P}(V_{x}^{G_{n}'}=A) \mathbb{P}(A \textrm{ is
directly connected to } b \textrm{ with a unique open edge}) \nonumber \\
&= \sum_{A} \mathbb{P}(V_{x}^{\Lambda_{l}}=A) n|A| (1-e^{-h/n})e^{-(|A|-1/n)h} \nonumber \\
&= n(e^{h/n}-1)\sum_{A} |A| \mathbb{P}(V_{x}^{\Lambda_{l}}=A)
e^{-|A|h} \nonumber \\ & = n(e^{h/n}-1) \frac{\partial
M_{x}^{\Lambda_{l}}}{\partial h}.
\end{align}
Here the sums are taken over all possible realizations $A$ of the
set of vertices $V_{x}^{G_{n}'}$. The second equality follows from
the independence of the bond variables in $G_{n}'$ and bond
variables which correspond to edges incident to $b$. The last
equality follows from (\ref{eq:
inequalities proof 2-1}).

The probability of the event $F_{2}$
is easily bounded from above using the BK-inequality
\begin{equation}\label{eq: inequalities proof 6}
\mathbb{P}(F_{2}) \leq \mathbb{P}(x \leftrightarrow b)^{2} =
(M_{x}^{\Lambda_{l}})^{2}.
\end{equation}
Now we bound the probability of the event $F_{3}$. Let
$C_{x}^{G_{n}' \backslash [y,z]_{i}}$ be the cluster in the graph
$G_{n}' \backslash [y,z]_{i}$, containing $x$ and $V_{x}^{G_{n}'
\backslash [y,z]_{i}}$ its vertex set. The event $F_{3}$ can be
partitioned with respect to realizations of
$V_{x}^{G_{n}'\backslash [y,z]_{i}}$. For a given set of vertices
$A$ we write $\left\{A \stackrel{d}{\leftrightarrow} b\right\}$
for the event that an edge between $b$ and some vertex in $A$ is
open, and $\left\{A \not \stackrel{d}{\leftrightarrow} b\right\}$
for the complement of this event. We obtain
$$
F_{3} = \bigcup_{(y,z) \in \Lambda_{l}^{2} \atop y \sim z}
\bigcup_{i=1}^{n} \bigcup_{A;y \in A}
\left\{V_{x}^{G_{n}'\backslash [y,z]_{i}}=A, [y,z]_{i} \textrm{ is
open}, A \not \stackrel{d}{\leftrightarrow} b, \left\{z
\leftrightarrow_{G_{n}\backslash A}b\right\} \circ \left\{z
\leftrightarrow_{G_{n}\backslash A}b\right\}\right\}.
$$
The union is taken over all possible realizations $A$ of the set
of vertices $V_{x}^{G_{n}'\backslash [y,z]_{i}}$. Here $G_{n}
\backslash A$ stands for the set of vertices of the graph $G_{n}$
which are not elements of $A$. The event
$\left\{V_{x}^{G_{n}'\backslash [y,z]_{i}}=A\right\}$ depends only
on the states of edges of $G_{n}'$ which have at least one
endpoint in $A$, but not on $[y,z]_{i}$. The event $\left\{A
\stackrel{d}{\leftrightarrow} b\right\}$ depends on the state of
edges between $b$ and the vertices in $A$. Finally, the event
$\left\{z \leftrightarrow_{G_{n}\backslash A}b\right\} \circ
\left\{z \leftrightarrow_{G_{n}\backslash A}b\right\}$ depends
only on the state of edges of $G_{n}$ which have no endpoints in
$A$. So we see that these events are independent and also
independent of the event $\left\{[y,z]_{i} \textrm{ is
open}\right\}$. Using this independence and the trivial fact
$\mathbb{P}([y,z]_{i} \textrm{ is open}) = (e^{\beta
/n}-1)\mathbb{P}([y,z]_{i} \textrm{ is closed})$, we get
\begin{multline}\label{eq: inequality proof 7}
\mathbb{P}(F_{3}) \leq (e^{\beta /n}-1) \sum_{(y,z) \in
\Lambda_{l}^{2} \atop y \sim z} \sum_{i=1}^{n} \sum_{A;y \in A}
\mathbb{P}(V_{x}^{G_{n}'\backslash
[y,z]_{i}}=A)\mathbb{P}([y,z]_{i} \textrm{ is closed}) \\
\mathbb{P}\big(\left\{z \leftrightarrow_{G_{n}\backslash
A}b\right\} \circ \left\{z \leftrightarrow_{G_{n}\backslash
A}b\right\}\big)\mathbb{P}(A \not \stackrel{d}{\leftrightarrow}
b).
\end{multline}
Notice that the BK inequality, Lemma \ref{lemma: cluster size
inequalities} a) and Lemma \ref{lemma: basic modification} imply
\begin{equation}\label{eq: inequality proof 8}
\mathbb{P}\big(\left\{z \leftrightarrow_{G_{n}\backslash
A}b\right\} \circ \left\{z \leftrightarrow_{G_{n}\backslash
A}b\right\}\big) \leq \mathbb{P}(z
\leftrightarrow_{G_{n}\backslash A}b)^{2} \leq \mathbb{P}(z
\leftrightarrow_{G_{n}\backslash
A}b)(M^{\Lambda_{l}}(\beta,h)+e^{-lh}).
\end{equation}
Inserting (\ref{eq: inequality proof 8}) into (\ref{eq: inequality
proof 7}) and using the independence again we obtain
\begin{multline}\label{eq: inequality proof 9}
\mathbb{P}(F_{3}) \leq (e^{\beta /n}-1)(M^{\Lambda_{l}}+e^{-lh})
\\
\sum_{(y,z)\in \Lambda_{l}^{2}} \sum_{i=1}^{n} \sum_{A;y \in A}
\mathbb{P}(V_{x}^{G_{n}'\backslash [y,z]_{i}}=A, [y,z]_{i}
\textrm{ is closed}) \mathbb{P}(A \not
\stackrel{d}{\leftrightarrow} b) \mathbb{P}(z
\leftrightarrow_{G_{n}\backslash A}b).
\end{multline}
In the last sum there are no contributions from $A$'s which
contain $z$, because, for such $A$'s, the set $\left\{z
\leftrightarrow_{G_{n}\backslash A}b\right\}$ is empty. So we can
take the sum over all possible realizations $A$ of $V_{x}^{G_{n}'
\backslash [y,z]_{i}}$ which contain $y$ but not $z$. For such
$A$'s it is clear that
$$
\left\{V_{x}^{G_{n}'\backslash [y,z]_{i}}=A, [y,z]_{i} \textrm{ is
closed}\right\} = \left\{V_{x}^{G_{n}'}=A\right\}.
$$
So (\ref{eq: inequality proof 9}) can be written in the form
$$
\mathbb{P}(F_{3}) \leq n(e^{\beta /n}-1)(M^{\Lambda_{l}}+e^{-lh})
\sum_{(y,z)\in \Lambda_{l}^{2} \atop y \sim z} \sum_{A;y \in
A}\mathbb{P}(V_{x}^{G_{n}'}=A) \mathbb{P}(A \not
\stackrel{d}{\leftrightarrow} b) \mathbb{P}(z
\leftrightarrow_{G_{n}\backslash A}b).
$$
The events $\left\{V_{x}^{G_{n}'}=A\right\}$, $\left\{A \not
\stackrel{d}{\leftrightarrow} b\right\}$ and $\left\{z
\leftrightarrow_{G_{n}\backslash A}b\right\}$ have the same
probabilities as the events $\left\{V_{x}^{\Lambda_{l}} =
A\right\}$, $\left\{A \cap B = \emptyset\right\}$ and $\left\{z
\leftrightarrow_{\Lambda_{l} \backslash A}B\right\}$, respectively.
Since the latter events are independent, we can write
\begin{equation*}
\mathbb{P}(F_{3}) \leq n(e^{\beta / n}-1)(M^{\Lambda_{l}} +
e^{-lh}) \sum_{(y,z) \in \Lambda_{l}^{2} \atop y \sim z} \sum_{A;
y \in A}\mathbb{P}(V_{x}^{\Lambda_{l}}=A, A\cap B = \emptyset, z
\leftrightarrow _{\Lambda_{l}\backslash A}B).
\end{equation*}
Using (\ref{eq: inequalities proof 1}), we obtain
\begin{equation}\label{eq: eq: inequality proof 10}
\mathbb{P}(F_{3}) \leq n(e^{\beta /n}-1)(M^{\Lambda_{l}}+e^{-lh})
\frac{\partial M_{x}^{\Lambda_{l}}}{\partial \beta}.
\end{equation}
Summing (\ref{eq: inequalities proof 5}), (\ref{eq: inequalities
proof 6}) and (\ref{eq: eq: inequality proof 10}) we get
$$
M_{x}^{\Lambda_{l}} \leq n(e^{h/n}-1) \frac{\partial
M_{x}^{\Lambda_{l}}}{\partial h} + (M_{x}^{\Lambda_{l}})^{2} +
n(e^{\beta /n}-1)(M^{\Lambda_{l}}+e^{-lh}) \frac{\partial
M_{x}^{\Lambda_{l}}}{\partial \beta}.
$$
For fixed $\beta$ and $h$ let $n$ go to $\infty$ and obtain
\begin{equation}\label{eq: inequality proof 11}
M_{x}^{\Lambda_{l}} \leq h \frac{\partial
M_{x}^{\Lambda_{l}}}{\partial h} + (M_{x}^{\Lambda_{l}})^{2} +
\beta (M^{\Lambda_{l}}+e^{-lh}) \frac{\partial
M_{x}^{\Lambda_{l}}}{\partial \beta}.
\end{equation}
Now sum (\ref{eq: inequality proof 11}) over $x \in \mathcal{F}$, use
the fact that $\displaystyle \sum_{x \in \mathcal{F}}
(M_{x}^{\Lambda_{l}})^{2} \leq (\sum_{x \in \mathcal{F}}
M_{x}^{\Lambda_{l}})^{2}$, and (\ref{eq: inequalities 2}) is
proved.
\end{proof}

\section{Site model}
\label{s-site-model}
In this part we shall explain how to obtain inequalities similar to those in (\ref{eq: inequalities 1}) and (\ref{eq: inequalities 2}), in the case of the site model. In the next section we shall use this inequalities to prove Theorem \ref{thm: main}. We will follow the arguments from the bond model and explain modifications necessary to proceed in the site case. In the following, just like before, $G$ will always denote a quasi-transitive graph. Throughout the section, $\overline{G'}$ will denote the subgraph induced by the set of vertices which lie in the subgraph $G'$ or have a neighbor in $G'$.

First we shall slightly change the notion of the cluster. We shall adopt the definitions from Section 7 in \cite{AizenmanB-87}.  Let $G$ be an arbitrary quasi-transitive graph, $x$ an arbitrary element of $G$ and $\omega$ an arbitrary site configuration. We define the modified cluster $\widetilde{C}_{x}(\omega)$ as the subgraph of $G$ induced by the following set of vertices:
$$
x \cup \left\{y; \textrm{there is a path from } x \textrm{ to } y \textrm{ including only open sites and } x\right\}.
$$
The distribution of the random variable $|\widetilde{C}_{x}|$ is clearly  different from the 
distribution of $|{C}_{x}|$. However, the probabilities $\mathbb{P}(|\widetilde{C}_{x}|=n), \, n \in \mathbb{N}\cup\left\{\infty\right\}$  will be proportional to $\mathbb{P}(|C_{x}|=n), \, n \in \mathbb{N}\cup\left\{\infty\right\}$. Therefore the critical values $p_{H}$ and $p_{T}$  will remain the same and the property of exponential decay below $p_{H}$ will be preserved. 

For any finite subgraph $G'$, the cluster $\widetilde{C}_{x}^{G'}(\omega)$ will be defined accordingly. Since the state of the site variable at some vertex $x$ is irrelevant for the properties of the cluster at $x$, we shall actually define the cluster $\widetilde{C}_{x}^{G'}(\omega)$ for all vertices $x$ and all subgraphs $G'$ such that $x \in \overline{G'}$.
Note that the relation ''$x$ lies in the cluster of $y$'' is not symmetric anymore and thus clusters can no longer be represented as connected components of some percolation subgraph.

The functions $M_{y}$, $M$ and their finite volume counterparts will be defined in the same way as before, where one just replaces $C_{x}$ by $\widetilde{C}_{x}$. Note that for these definitions as well as for Lemma \ref{lemma: OP definition}, and Propositions \ref{prop: OP properties} and \ref{prop: probability limits} one only needs to define the notion of cluster, while the underlying model is completely irrelevant. This is why these results transfer directly to the site percolation setting. Lemma \ref{lemma: cluster size inequalities} holds in the site model as well and the proof remains practically the same. This also holds for Lemma \ref{lemma: basic modification} and Proposition \ref{lemma: convergence}. The notion of blue sites is introduced in the same way as in section \ref{s-order-parameter}. Similarly as above, the state of the site variable at some vertex $x$ is irrelevant for the connectedness to any other vertex and so by $x \leftrightarrow _{G'} B$ we will denote the event that $x$ is connected to some blue site by a path in which all vertices, except maybe $x$, lie in the subgraph $G'$. Proposition \ref{prop: blue sites} remains unchanged in the site setting.

Now we establish differential inequalities for the site model. 
In the proofs one has to be careful not to mix two types of site variables, those which correspond to the percolation process and those which correspond to blue sites. In particular, pivotality will refer only to percolation variables. Similarly as before $\Lambda_{l}\backslash G'$ will denote the subgraph of $\Lambda_{l}$ obtained by deleting all vertices in $G'$ and all edges which are incident to some vertex in $G'$.

The formula in Proposition \ref{prop: inequality 1} remains the same in the site percolation setting.

\begin{proposition}\label{prop: inequality 1 site}
There exists a constant $K>0$ such that 
\begin{equation}\label{eq: inequalities 1 site}
\frac{\partial M^{\Lambda_{l}}}{\partial \beta} \leq K
(M^{\Lambda_{l}}+e^{-lh})\frac{\partial M^{\Lambda_{l}}}{\partial
h}, \textrm{ for all positive } \beta \textrm{ and } h.
\end{equation}
\end{proposition}

\begin{proof}
The proof is analogous to the proof of Proposition \ref{prop: inequality 1}. Using the Russo formula one obtains 
$$
\frac{\partial M_{x}^{\Lambda_{l}}}{\partial \beta} = \sum_{y
\in \Lambda_{l}} \mathbb{P}(y \textrm{ is closed},y
\textrm{ is pivotal for the event } \left\{ x \leftrightarrow
_{\Lambda_{l}}B \right\}).
$$
We notice that
\begin{multline*}
\left\{y \textrm{ is closed},y \textrm{ is pivotal for the
event } \left\{ x \leftrightarrow _{\Lambda_{l}}B
\right\}\right\}\\
= \left\{\widetilde{C}_{x}^{\Lambda_{l}} \cap B = \emptyset, y \in
\overline{\widetilde{C}_{x}^{\Lambda_{l}}}, y \leftrightarrow _{\Lambda_{l}\backslash
\overline{\widetilde{C}_{x}^{\Lambda_{l}}}}B\right\}.
\end{multline*}
Now using independence and presenting the above events as unions over possible realizations of $\widetilde{C}_{x}^{\Lambda_{l}}$ we obtain by means of the site version of Proposition \ref{prop: blue sites}
\begin{equation}\label{eq: inequalities 1-2 site}
\frac{\partial M_{x}^{\Lambda_{l}}}{\partial \beta} = \sum_{y
\in \Lambda_{l}} \sum_{{A \subset \Lambda_{l}} \atop {y \in
\overline{A} \backslash A}}\mathbb{P}(\widetilde{C}_{x}^{\Lambda_{l}}=A) \, M_{y}^{\Lambda_{l}\backslash \overline{A}}
\, e^{-h|A|}.
\end{equation}
Using the site versions of Lemmata \ref{lemma: cluster size inequalities} a) and \ref{lemma: basic modification} we get
\begin{align}\label{eq: inequalities proof 2 site}
\frac{\partial M_{x}^{\Lambda_{l}}}{\partial \beta} & \leq (M^{\Lambda_{l}}+e^{-lh})\sum_{A \subset \Lambda_{l}}|\overline{A} \backslash A|\mathbb{P}(\widetilde{C}_{x}^{\Lambda_{l}}=A) e^{-h|A|} \nonumber \\ & \leq 
K(M^{\Lambda_{l}}+e^{-lh}) \sum_{A \subset \Lambda_{l}}|A| \,
\mathbb{P}(\widetilde{C}_{x}^{\Lambda_{l}}=A) \, e^{-h|A|} \nonumber \\ & = K(M^{\Lambda_{l}}+e^{-lh}) \frac{\partial M_{x}^{\Lambda_{l}}}{\partial h},
\end{align}
where $K$ is again the maximal vertex degree in $G$. The last equality in (\ref{eq: inequalities proof 2 site}) is proven just like its bond analogue (see (\ref{eq: inequalities proof 2})).
\end{proof}

In Proposition \ref{prop: inequality 2} we will modify one term on the left hand side. This is due to the fact that we will not change the graph $\Lambda_{l}$ as dramatically as in the proof of Proposition \ref{prop: inequality 2}.
\begin{proposition}\label{prop: inequality 2 site}
The finite volume order parameter satisfies the following inequality
\begin{equation}\label{eq: inequalities 2 site}
M^{\Lambda_{l}} \leq h \frac{\partial M^{\Lambda_{l}}}{\partial h}
+ \big(M^{\Lambda_{l}}\big)^{2} + (e^{\beta}-1)
(M^{\Lambda_{l}}+e^{-lh})\frac{\partial M^{\Lambda_{l}}}{\partial
\beta}, \textrm{ for all positive } \beta \textrm{ and } h.
\end{equation}
\end{proposition}

\begin{proof}
As in the proof of Proposition \ref{prop: inequality 2} we will change both our graph and probability space. Rather than considering the probability space of blue sites we shall add a new vertex $b$ to the graph $\Lambda_{l}$. This vertex will be connected to each vertex in $\Lambda_{l}$ by exactly $n$ edges. Notice that the new graph, which will be denoted by $G_{n}$, contains $\Lambda_{l}$ as a subgraph. We shall now consider a mixed site-bond percolation model in which each vertex of $\Lambda_{l}$ is open with probability $1-e^{-\beta}$ and in which each bond incident to $b$ is open with probability $1-e^{-h/n}$. Again we shall compare this process with the usual percolation model with blue sites which are generated on each vertex with probability $1-e^{-h}$. For the relation between these two processes see the discussion at the beginning of the proof of Proposition \ref{prop: inequality 2}. The event $\left\{x \leftrightarrow_{G_{n}} b\right\}$ can be partitioned into following disjoint events:
\begin{align*}
&F_{1}:=\left\{\textrm{There is a unique open edge which connects some vertex of } \widetilde{C}_{x}^{\Lambda_{l}} \textrm{
with } b\right\},
\\ &
F_{2}:=\left\{x \leftrightarrow _{G_{n}} b\right\} \circ \left\{x
\leftrightarrow _{G_{n}}b\right\}\\& \ \ \ \ = \left\{\textrm{There are two paths from } x \textrm{ to } b, \textrm{ which have no common vertices other than } x \textrm{ and } b\right\},
\\ &
F_{3}:=\bigcup_{y \in \Lambda_{l}} \left\{y \in \widetilde{C}_{x}^{\Lambda_{l}} \textrm{ is open and pivotal for
} \left\{x \leftrightarrow _{G_{n}}b\right\}, \ \left\{y
\leftrightarrow _{G_{n}} b\right\} \circ
\left\{y \leftrightarrow _{G_{n}}b\right\}\right\},
\end{align*}
The probability of the first event can be calculated in the same way as in the bond case. We obtain:
\begin{equation}\label{eq: inequalities proof 5 site}
\mathbb{P}(F_{1}) = n(e^{h/n}-1)\frac{\partial M_{x}^{\Lambda_{l}}}{\partial h}.
\end{equation}
Since $\mathbb{P}(x \leftrightarrow _{G_{n}} b)=\mathbb{P}(x \leftrightarrow _{\Lambda_{l}} B)$, the BK inequality implies
\begin{equation}\label{eq: inequalities proof 6 site}
\mathbb{P}(F_{2}) \leq \mathbb{P}(x \leftrightarrow _{G_{n}}b)^{2} = (M_{x}^{\Lambda_{l}})^{2}.
\end{equation}
The last event can be rewritten as 
$$
F_{3} = \bigcup_{y \in \Lambda_{l}} \bigcup_{A; y \in \overline{A} \backslash A} \left\{\widetilde{C}_{x}^{\Lambda_{l}\backslash \left\{y\right\}} = A, A \not \stackrel{d}{\leftrightarrow} b, y \textrm{ open }, \left\{y \leftrightarrow _{G_{n} \backslash \overline{A}} b\right\} \circ \left\{y \leftrightarrow _{G_{n} \backslash \overline{A}} b\right\}\right\},
$$
where the second union is taken over all possible realizations $A$ of $\widetilde{C}_{x}^{\Lambda_{l}\backslash \left\{y\right\}}$, for which $y \in \overline{A} \backslash A$ (note that $\overline{A}$ is defined as the closure of $A$ in $\Lambda_{l}$).
Considering the above formula and using independence we obtain
\begin{align}
\label{eq: inequality proof 9 site}
\mathbb{P}(F_{3}) 
& \leq (e^{\beta}-1) \sum_{y \in \Lambda_{l}} \sum_{A;y \in \overline{A} \backslash A} \mathbb{P}(\widetilde{C}_{x}^{\Lambda_{l}\backslash \left\{y\right\}}=A, y \textrm{ is closed}) 
\\ \nonumber
& \hspace*{10em} 	\mathbb{P}\Big(\left\{y \leftrightarrow _{G_{n} \backslash \overline{A}} b\right\} \circ \left\{y \leftrightarrow _{G_{n} \backslash \overline{A}} b\right\}\Big) \mathbb{P}(A \not \stackrel{d}{\leftrightarrow} b) 
\\ \nonumber
& \leq (e^{\beta}-1)(M^{\Lambda_{l}}+e^{-lh}) \sum_{y \in \Lambda_{l}}\sum_{A;y \in \overline{A}\backslash A} \mathbb{P}(\widetilde{C}_{x}^{\Lambda_{l}}=A)M_{y}^{\Lambda_{l} \backslash \overline{A}}e^{-h|A|}\\ \nonumber
&= (e^{\beta}-1)(M^{\Lambda_{l}}+e^{-lh}) \frac{\partial M_{x}^{\Lambda_{l}}}{\partial \beta}
\end{align}
In the second inequality we used the BK inequality and the site versions of Lemma \ref{lemma: cluster size inequalities} a), Lemma \ref{lemma: basic modification} and Proposition \ref{prop: blue sites} in the same way as in the proof of Proposition \ref{prop: inequality 2}. In the last equality we used (\ref{eq: inequalities 1-2 site}). Now the result follows after taking the sum of (\ref{eq: inequalities proof 5 site}), (\ref{eq: inequalities proof 6 site}) and (\ref{eq: inequality proof 9 site}) and letting $n$ tend to $\infty$.

\end{proof}

\section{Completion of the proof of Theorem \ref{thm: main}}
\label{s-completion}

In this section we will complete the proof of our main result,
using the differential inequalities (\ref{eq: inequalities 1}) and
(\ref{eq: inequalities 2}).

The next result will be useful in  the proof of Lemma \ref{lemma:
critical asymptotics}. It is a special case of Lemma 4.1 in
\cite{AizenmanB-87}.

\begin{lemma}\label{lemma: from AB}
Let $M \colon \mathbb{R}^{+} \to \mathbb{R}$ be an increasing
differentiable function of $h$ obeying
$$
\lim_{h \downarrow 0}M(h) = 0, \ \ \ \lim_{h \downarrow 0}
\frac{M(h)}{h} = \infty \ \ \ and \ \ \ M \leq h\frac{dM}{dh} +
M^{2} + kM^{2}\frac{dM}{dh}, \textrm{ for all } h>0,
$$
for some positive constant $k$. Then there exists a constant $c>0$
such that for all $h>0$ small enough we have
$$
M(h) \geq c \sqrt{h}.
$$
\end{lemma}

Lemma \ref{lemma: critical asymptotics} and Proposition \ref{prop:
final proposition} are the final steps of the proof of Theorem
\ref{thm: main}. They correspond to Theorem 4.2 and Lemma 5.1 from
\cite{AizenmanB-87}. We just need some adjustments in the proof of
Proposition \ref{prop: final proposition} because of somewhat
different differential inequalities. As an aside we obtain an
upper bound on the critical exponent defined as
\begin{equation} \label{eqn: critical exponent}
\delta :=\liminf_{h \downarrow 0} \frac{\ln h}{M(\beta_T,h)}.
\end{equation}

\begin{lemma}\label{lemma: critical asymptotics}
There is a constant $c > 0$, such that for $h>0$ small enough
\begin{equation}\label{eq: critical asymptotics 0}
M(\beta_{T},h) \geq c \sqrt{h}.
\end{equation}
In particular, the critical exponent \eqref{eqn: critical
exponent} obeys $\delta \ge 2$.
\end{lemma}

\begin{proof}
First notice that $M$ satisfies the following differential
inequality
\begin{equation}\label{eq: critical asymptotics 1}
M(\beta,h) \leq h\frac{\partial M}{\partial h}(\beta,h) +
M^{2}(\beta,h) + K\beta M^{2}(\beta,h) \frac{\partial M}{\partial
h} (\beta,h),
\end{equation}
in the bond case, respectively
\begin{equation}\label{eq: critical asymptotics 1'}
M(\beta,h) \leq h\frac{\partial M}{\partial h}(\beta,h) +
M^{2}(\beta,h) + K(e^{\beta}-1) M^{2}(\beta,h) \frac{\partial M}{\partial
h} (\beta,h),
\end{equation}
in the site case. These inequalities can be proven by inserting (\ref{eq: inequalities 1}) into
(\ref{eq: inequalities 2}) in the bond case (or inserting (\ref{eq: inequalities 1 site}) into (\ref{eq: inequalities 2 site}) in the site case) and letting $l$ go to $\infty$. If
$\lim_{h \downarrow 0}M(\beta_{T},h)
> 0$ there is nothing to prove. Suppose $\lim_{h \downarrow 0}M(\beta_{T},h)
= 0$. This implies by Proposition \ref{prop: probability limits} that
$\mathbb{P}_{\beta_{T}}(|C_{x}|=\infty)=0$, for all vertices $x$ and
$$
\lim_{h \downarrow 0} \frac{\partial M}{\partial h}(\beta_{T},h) =
\sum_{x \in \mathcal{F}}
\mathbb{E}_{\beta_{T}}(|C_{x}|;|C_{x}|<\infty) = \sum_{x \in
\mathcal{F}} \mathbb{E}_{\beta_{T}}(|C_{x}|) = \infty,
$$
by Proposition \ref{prop: mean cluster size}. Now the Mean Value Theorem implies
\begin{equation}\label{eq: critical asymptotics 2}
\lim_{h \downarrow 0}\frac{M(\beta_{T},h)}{h} = \infty.
\end{equation}
In view of (\ref{eq: critical asymptotics 1}) (respectively (\ref{eq: critical asymptotics 1'}) in the site case) and (\ref{eq: critical asymptotics 2}),
the claim follows directly from Lemma \ref{lemma: from AB}.
\end{proof}

Except for having to control the term $e^{-hl}$ in (\ref{eq: inequalities 1 site}) and (\ref{eq: inequalities 2 site}), the proof of the next proposition is the same as the proof of Lemma 5.1 in \cite{AizenmanB-87}.

\begin{proposition}\label{prop: final proposition}
For any $\beta' > \beta_{T}$ we can find a positive constant $d >
0$ such that
\begin{equation}\label{eq: final proposition 1}
\lim_{h \downarrow 0} M(\beta,h) \geq d(\beta-\beta_{T})
\end{equation}
holds for every $\beta \in [\beta_{T},\beta']$.
\end{proposition}

\begin{proof}
Let's consider the bond case first. Change the variables $(\beta,h)$ to $(\beta,\ln h)$, i.e.~define
$u:=\ln h$ and $\widetilde{M}^{\Lambda_{l}}(\beta,u) =
M^{\Lambda_{l}}(\beta,h)$. Now $\displaystyle \frac{\partial \ln
\widetilde{M}^{\Lambda_{l}}}{\partial u}(\beta,u) =
\frac{h}{M^{\Lambda_{l}}(\beta,h)}\frac{\partial
M^{\Lambda_{l}}}{\partial h}(\beta,h)$ and so (\ref{eq:
inequalities 2}) can be rewritten as
\begin{equation}\label{eq: final proposition proof 1}
1 \leq \frac{\partial \ln \widetilde{M}^{\Lambda_{l}}}{\partial
u}(\beta,u) + \widetilde{M}^{\Lambda_{l}}(\beta,u) + \beta \Big(1
+
\frac{e^{-le^{u}}}{\widetilde{M}^{\Lambda_{l}}(\beta,u)}\Big)\frac{\partial
\widetilde{M}^{\Lambda_{l}}}{\partial \beta}(\beta,u),
\end{equation}
for every $u \in \mathbb{R}$. Now fix some $0 < h_{1} < h_{2}$,
define $u_{1}:=\ln h_{1}$ and $u_{2}:=\ln h_{2}$ and integrate
(\ref{eq: final proposition proof 1}) over the rectangle
$[\beta_{T},\beta_{1}]\times[u_{1},u_{2}]$, where $\beta_{1}$ is
an arbitrary real number between $\beta_{T}$ and $\beta'$. Using
the fact that $\widetilde{M}^{\Lambda_{l}}$ is increasing in both
$\beta$ and $u$ and switching back to $h$, we get
\begin{multline*}
(\beta_{1}-\beta_{T})\ln \frac{h_{2}}{h_{1}} \leq
(\beta_{1}-\beta_{T})\ln \frac{M^{\Lambda_{l}}(\beta_{1},h_{2})}
{M^{\Lambda_{l}}(\beta_{T},h_{1})} + (\beta_{1}-\beta_{T})\ln
\frac{h_{2}}{h_{1}} M^{\Lambda_{l}}(\beta_{1},h_{2}) \\ + \beta'
\ln \frac{h_{2}}{h_{1}}
\Big(1+\frac{e^{-lh_{1}}}{M^{\Lambda_{l}}(\beta_{T},h_{1})}\Big)(M^{\Lambda_{l}}(\beta_{1},h_{2})
- M^{\Lambda_{l}}(\beta_{T},h_{1})).
\end{multline*}
Let $l$ go to $\infty$ and obtain
\begin{multline}\label{eq: final proposition 2}
\beta_{1}-\beta_{T} \leq (\beta_{1}-\beta_{T}) \frac{\ln \frac{
M(\beta_{1},h_{2})}{M(\beta_{T},h_{1})}}{\ln \frac{h_{2}}{h_{1}}}
+ (\beta_{1}-\beta_{T}) M(\beta_{1},h_{2}) + \beta'
(M(\beta_{1},h_{2}) - M(\beta_{T},h_{1})).
\end{multline}
Now notice
\begin{equation}\label{eq: final proposition 3}
\frac{\ln \frac{M(\beta_{1},h_{2})}{M(\beta_{T},h_{1})}}{\ln
\frac{h_{2}}{h_{1}}} = \frac{\ln M(\beta_{1},h_{2}) - \ln
M(\beta_{T},h_{1})}{\ln h_{2} - \ln h_{1}} = \frac{\frac{\ln
M(\beta_{1},h_{2})}{\ln h_{1}} - \frac{\ln M(\beta_{T},h_{1})}{\ln
h_{1}}}{\frac{\ln h_{2}}{\ln h_{1}} - 1}.
\end{equation}
Using Lemma \ref{lemma: critical asymptotics} and (\ref{eq: final
proposition 3}) we get
\begin{equation}\label{eq: final proposition 4}
\limsup_{h_{1} \downarrow 0} \frac{\ln
\frac{M(\beta_{1},h_{2})}{M(\beta_{T},h_{1})}}{\ln
\frac{h_{2}}{h_{1}}} \leq \frac{1}{2}.
\end{equation}
Inserting (\ref{eq: final proposition 4}) to (\ref{eq: final
proposition 2}) and letting $h_{1} \downarrow 0$ leads to
$$
\frac{1}{2} (\beta_{1}-\beta_{T}) \leq
M(\beta_{1},h_{2})(\beta_{1}-\beta_{T}+ \beta') - \beta'
\lim_{h_{1} \downarrow 0}M(\beta_{T},h_{1}) \leq
(2\beta'-\beta_{T})M(\beta_1,h_{2}).
$$
Let $h_{2} \downarrow 0$ and the proof is over.

In the site model we start by changing the variables $(\beta, h)$ to $(p,u):=(1-e^{-\beta}, \ln h)$.
In other words, this time  we define $\widetilde{M} ^{\Lambda_{l}}\colon ]0,1[ \times \mathbb{R} \to \mathbb{R}$ such that  $\widetilde{M}^{\Lambda_{l}}(p,u):=M^{\Lambda_{l}}(\beta,h)$. Now (\ref{eq: inequalities 2 site}) can be rewritten as
\begin{equation}\label{eq: final proposition proof 1 site}
1 \leq \frac{\partial \ln \widetilde{M}^{\Lambda_{l}}}{\partial
u}(p,u) + \widetilde{M}^{\Lambda_{l}}(p,u) + p \Big(1
+
\frac{e^{-le^{u}}}{\widetilde{M}^{\Lambda_{l}}(p,u)}\Big)\frac{\partial
\widetilde{M}^{\Lambda_{l}}}{\partial p}(p,u).
\end{equation}
This inequality replaces (\ref{eq: final proposition proof 1}) but has the same form. Thus the proof continues the same way as in the bond case after making the transformations $\beta_{T} \mapsto p_{T}:=1-e^{-\beta_{T}}$ and $\beta' \mapsto p':=1-e^{-\beta'}$.
\end{proof}

\begin{proof}[Proof of Theorem \ref{thm: main}]
Proposition \ref{prop: final proposition} tells us that $\lim_{h
\downarrow 0}M(\beta,h)$ is positive as soon as $\beta >
\beta_{T}$. In the view of Proposition \ref{prop: probability
limits} this proves the main result.
\end{proof}

\section{Extension of results to general bond models}

In this section we will explain how the methods presented above can be
applied to more general bond percolation models on
quasi-transitive graphs. The model we present here is the
\emph{partially oriented long-range model} which was considered in
\cite{AizenmanB-87} for the lattice case.

Assume that $G=(V,E)$ is again a quasi-transitive graph, with some
fixed fundamental domain $\mathcal{F}$. Now make the graph
complete, that is connect each pair of vertices
$\left\{x,y\right\}$ with an unoriented edge $[x,y]$. Moreover,
connect $x$ and $y$ with two oriented edges $\left[x,y\right>$
(oriented from $x$ to $y$) and $\left[y,x\right>$ (oriented from
$y$ to $x$). The distance function on the vertices is the one inherited from
the graph $G$. Thus it makes sense to define the length of an edge
as the distance (in $G$) between its endvertices. Paths in our graph can
contain both oriented and unoriented edges, but the orientation of
oriented edges must be in accordance with the orientation of the considered
path.

On the complete graph the usual nearest neighbor bond percolation is
uninteresting, because any parameter $p>0$ will correspond to the
supercritical phase. To avoid this triviality, one has to
introduce certain damping of the probabilities that $x$ and $y$
are connected, as the distance between $x$ and $y$ goes to
infinity. This is done by introducing for each pair of vertices $(x,y)$ two positive parameters
$J_{[x,y]}$ and $J_{\left[x,y\right>}$. The unoriented edge $[x,y]$
will be open with probability $1-e^{-\beta J_{[x,y]}}$ and
the oriented edge $\left[x,y\right>$ will be open with probability
$1-e^{-\beta J_{\left[x,y\right>}}$. Of course, we assume that all
these events are mutually independent and thus the product probability space can be
constructed similarly as before. The structure of the quasi-transitive
graph $G$ is reflected through the invariance of the parameters
$J$: we assume that the parameters $J$ are invariant under the
automorphisms of the graph $G$. In other words, $J_{[\gamma
x,\gamma y]}= J_{[x,y]}$ and $J_{\left[\gamma x,\gamma y\right>} =
J_{\left[x,y\right>}$, for all $\gamma \in Aut(G)$ and all
vertices $x$ and $y$. Next we define $J_{x}:=\sum_{y \in
V}(J_{[x,y]}+J_{\left[x,y\right>})$. To avoid the triviality
mentioned above, we will assume that
\begin{equation}\label{eq: long model 1}
J_0:=\sup_{x \in V} J_{x} =\max_{x \in \mathcal{F}} J_{x} <
\infty.
\end{equation}
Without this assumption, some vertices would be directly connected
with infinitely many other vertices almost surely.

The subgraphs $\Lambda_{l}$ are also defined similarly as before, using the distance function of the original graph $G$. The
vertex set remains unchanged, but for the set of edges we take all
possible oriented and unoriented edges between pairs of vertices contained in $\Lambda_l$.
Percolation on the graph $\Lambda_{l}$ inherits the probabilities
for edges to be open from the percolation on the whole graph.

Since the graph contains oriented edges, the relation "being connected
in a percolation subgraph" defined on the set of vertices is not
symmetric any more and thus the notion of the connected components
is now meaningless. However, the percolation cluster containing
some vertex $x$ can be defined in a natural way, as the graph
$C_{x}(\omega)$ for which the vertex set is the set of all vertices
which can be reached from $x$ by an open path. The edge set is defined
as the set of all open edges between vertices of $C_{x}(\omega)$.
A percolation cluster $C_{x}^{\Lambda_{l}}(\omega)$ in $\Lambda_l$ is defined
similarly. Using this new definition of clusters, the order parameter $M$
and the finite volume order parameter $M^{\Lambda_{l}}$ can be defined in
the same way as before. The probabilistic interpretation with colored
sites is also applicable just as before, since Proposition
\ref{prop: blue sites} is true in this setting, too.
\medskip

The critical parameters $\beta_{T}$ and $\beta_{H}$ are defined in the
same way as before. In the nearest neighbor model, the fact that
these values are well defined relied on the Fundamental Tools
presented in Section \ref{s-basics}. Both Definition \ref{def: events} and the
Fundamental Tools can be generalized to the present model in a
natural way. These results are scattered in the literature. For
example, a general Russo inequality can be found in \cite{Russo-81}
and a general BK inequality can be found in \cite{vandenBergK-85}.
For more explanations, one can
also look at the arguments in \cite{AizenmanB-87} regarding this
general model. Just as before, these generalizations imply the
fact that the critical parameters are well defined. The generalizations of the
Fundamental Tools also ensure that the basic properties of the order
parameter remain valid in the new model. Lemma \ref{lemma: OP
definition} and Propositions \ref{prop: OP properties} and
\ref{prop: probability limits} are still valid in the new setting. One
can easily convince oneself that this is also true for parts a)
and b) of Lemma \ref{lemma: cluster size inequalities}. However,
we need to be more careful with part c).
Rather than the equality stated in this part, in the general model we obtain an inequality formulated in Lemma \ref{lemma: long model 1} below.
This inequality will be used
in Lemma \ref{lemma: long range 2} which replaces Lemma
\ref{lemma: basic modification}. Under the same assumptions as in
Lemma \ref{lemma: basic modification}, Lemma \ref{lemma: long
range 2} gives the following inequality
\begin{equation}\label{eq: long model 2}
M_{y}^{\Lambda_{l}}(\beta,h) \leq M_{x}^{\Lambda_{l}}(\beta,h) +
f_{l}(\beta,h),
\end{equation}
where $(f_{l})_{l\in \NN}$ is some sequence of positive continuous
functions which converges to zero locally uniformly for $l\to\infty$.
Notice that this bound is sufficient to prove differential
inequalities similar to those in (\ref{eq: inequalities 1}) and
(\ref{eq: inequalities 2}). Namely, one obtains the following
differential inequalities
\begin{align}\label{eq: long model 3}
\frac{\partial M^{\Lambda_{l}}}{\partial \beta} &\leq K
(M^{\Lambda_{l}}+f_{l})\frac{\partial M^{\Lambda_{l}}}{\partial h}, \quad \text{ and }
\\
\label{eq: long model 4}
M^{\Lambda_{l}} &\leq h \frac{\partial M^{\Lambda_{l}}}{\partial h}
+ \big(M^{\Lambda_{l}}\big)^{2} + \beta
(M^{\Lambda_{l}}+f_{l})\frac{\partial M^{\Lambda_{l}}}{\partial
\beta}.
\end{align}
These inequalities are sufficient to conclude the equality of the
critical values $\beta_{H}=\beta_{T}$. More precisely, the proof of Lemma
\ref{lemma: critical asymptotics} extends to our new setting.
For this we
use the fact that Proposition \ref{prop: mean cluster size} is also
true in this general model (this also follows from the proof of
Lemma 3.1 in \cite{AizenmanN-84}). Proposition \ref{prop: final
proposition} still gives the main result, since its proof does not
require any special form of the functions $(f_{l})$, but only the
fact, that they decay locally uniformly for $l\to\infty$.
\medskip

Now we state the mentioned inequality which replaces part c) of Lemma
\ref{lemma: cluster size inequalities}.
\begin{lemma}\label{lemma: long model 1}
There exists a nondecreasing sequence of positive integers
$(n_{l})_{l\in\NN}$, which converges to infinity and a sequence of
positive continuous functions $(g_{l})_{l\in\NN}$, $g_{l} \colon
]0,\infty[ \to \mathbb{R}$ which converges to zero
locally uniformly for $l \to \infty$, such that the following inequality holds for
any $x \in \mathcal{F}$, any $l\in\NN$, and any positive integer $1 \leq k \leq
n_{l}$
\begin{equation}\label{eq: long model 5}
\mathbb{P}(|C_{x}|\geq k) \leq \mathbb{P}(|C_{x}^{\Lambda_{l}}|
\geq k) + g_{l}(\beta), \textrm{ for all } \beta \in ]0,\infty[.
\end{equation}
\end{lemma}

\begin{proof}
We follow the arguments in the proof of Lemma A.3 from
\cite{AizenmanB-87}. For any positive real $r$ define
$$
J_{r}:=\max_{x \in \mathcal{F}}\sum_{y \in V \atop d(x,y)\geq
r}(J_{[x,y]} + J_{\left[x,y\right>}).
$$
Since $J_0$ is finite, $\lim_{r\to\infty} J_r=0$. In the following
we will use the identification from Remark \ref{rem: 1} in our new
setting. We have to estimate $\mathbb{P}(|C_{x}|\geq k,
|C_{x}^{\Lambda_{l}}|<k)$. For any $\omega \in \left\{|C_{x}|\geq
k, |C_{x}^{\Lambda_{l}}|<k\right\}$ there exists a path consisting
of at most $k$ edges which connects some $x$ with some vertex
outside $\Lambda_{l}$. This path connects two vertices which are
at distance greater or equal to $l$. Thus there has to be an edge
in this path which has length greater or equal to $l/k$. To reach
this edge we have to make $j$ steps in the path, for some $j$ such
that $0 \leq j \leq k-1$. The probability that there exists an
open edge of length greater or equal than $l/k$ which can be
reached from $x$ by an open open path of length $j$, can be bounded
above by $\beta J_{l/k}\, (\beta J_0)^{j} $. Here we used the
inequality $1-e^{-t} \leq t$ for any positive $t$. So from the
arguments above we deduce that
\begin{equation}\label{eq: long model 5-1}
\mathbb{P}(|C_{x}|\geq k, |C_{x}^{\Lambda_{l}}|<k) \leq \beta
J_{l/k}\sum_{j=1}^{k}(\beta J_0)^{j-1}.
\end{equation}
For any positive integer $n$ define $\displaystyle
K_{n}:=n\sum_{k=1}^{n}(nJ_0)^{k-1}$. Since $\lim_{r \to
\infty}J_{r}=0$, we can find an increasing sequence of positive
integers $L_{n}$ such that
\begin{equation}\label{eq: long model 6}
\lim_{n \to \infty}J_{L_{n}/n}K_{n} = 0.
\end{equation}
Now define $n_{l}:=\max \left\{n; L_{n} \leq l\right\}$ and
$\displaystyle g_{l}(\beta):=\beta J_{l/n_{l}}
\sum_{k=1}^{n_{l}}(\beta J_0)^{k-1}$. From (\ref{eq: long model
6}) it is clear that $\lim_{l \to \infty}g_{l}(\beta) = 0$ locally
uniformly on $\mathbb{R}^{+}$. The other claimed properties of the
sequences $(n_{l})_{l}$ and $(g_{l})_{l}$ are obvious. Using
(\ref{eq: long model 5-1}) and the fact that $r \mapsto J_{r}$ is
a non-increasing function we obtain $\mathbb{P}(|C_{x}|\geq k,
|C_{x}^{\Lambda_{l}}|<k) \leq g_{l}(\beta)$. This proves the
lemma.
\end{proof}
Using the previous result and part a) of Lemma 12 which, as we
said, still holds in our new setting, one easily obtains $\lim_{l \to
\infty}\mathbb{P}(|C_{x}^{\Lambda_{l}}| \geq k) =
\mathbb{P}(|C_{x}|\geq k)$. This can be used to prove the
pointwise convergence of the finite volume order parameter to the
order parameter and the same claim for the partial derivative in $h$,
that is Proposition \ref{lemma: convergence}.
\smallskip

Using Lemma \ref{lemma: long model 1} one can easily obtain an
inequality as in (\ref{eq: long model 2}).

\begin{lemma}\label{lemma: long range 2}
Let $y$ be a vertex of $\Lambda_{l}$ and $x$ be the unique element
of $\mathcal{F}$ in the same orbit as $y$. Then there exists a
sequence of positive continuous functions $(f_{l})_{l\in\NN}$,
converging locally uniformly to $0$ for $l\to\infty$, such that the following
inequality holds
\begin{equation}\label{eq: long model 7}
M_{y}^{\Lambda_{l}}(\beta,h) \leq M_{x}^{\Lambda_{l}}(\beta,h) +
f_{l}(\beta,h), \textrm{ for all } (\beta,h) \in ]0,\infty[^{2}.
\end{equation}
\end{lemma}

\begin{proof}
Using a similar decomposition as in the proof of Lemma \ref{lemma:
basic modification}, and then Lemma \ref{lemma: long model 1} and Lemma
\ref{lemma: OP definition} we get
\begin{align}\label{eq: long model 8}
M_{y}^{\Lambda_{l}}(\beta,h) & =
\sum_{k=1}^{n_{l}}\mathbb{P}(|C_{y}^{\Lambda_{l}}|\geq
k)(e^{-(k-1)h}-e^{-kh}) + \sum_{n_{l}+1 \leq k <
\infty}\mathbb{P}(|C_{y}^{\Lambda_{l}}|\geq
k)(e^{-(k-1)h}-e^{-kh}) \nonumber
\\ & \leq \sum_{k=1}^{n_{l}}\mathbb{P}(|C_{x}^{\Lambda_{l}}|\geq
k)(e^{-(k-1)h}-e^{-kh}) +
\sum_{k=1}^{n_{l}}g_{l}(\beta)(e^{-(k-1)h}-e^{-kh}) + e^{-n_{l}h}
\nonumber
\\ & \leq M_{x}^{\Lambda_{l}}(\beta,h) + g_{l}(\beta) +
e^{-n_{l}h}.
\end{align}
Since $(n_{l})_{l}$ converges to infinity, the claim of the lemma is proven.
\end{proof}

\def\cprime{$'$}\def\polhk#1{\setbox0=\hbox{#1}{\ooalign{\hidewidth
  \lower1.5ex\hbox{`}\hidewidth\crcr\unhbox0}}}

\end{document}